\documentclass[a4paper,fleqn]{cas-sc}

\usepackage{setspace}
\usepackage{bm}
\usepackage{multirow}
\usepackage[numbers,sort&compress]{natbib}
\usepackage{graphicx}
\usepackage{float}
\usepackage{geometry}
\usepackage{enumitem}
\usepackage{amsmath}
\usepackage{times}
\usepackage{caption}

\usepackage[ruled]{algorithm2e}
\usepackage[title]{appendix}%

\usepackage{multirow}

\usepackage{lineno}
\def\tsc#1{\csdef{#1}{\textsc{\lowercase{#1}}\xspace}}
\tsc{WGM}
\tsc{QE}
\tsc{EP}
\tsc{PMS}
\tsc{BEC}
\tsc{DE}

\begin{document}
\onehalfspacing
\let\WriteBookmarks\relax
\def\floatpagepagefraction{1}
\def\textpagefraction{.001}
\shorttitle{Diffusion Model-Augmented Construction of Three-Body Earth-Moon Transfers}
\shortauthors{Shuyue Fu et~al.}

\title [mode = title]{Three-Body Earth-Moon Transfers with Different Departure/Arrival Orbital Altitudes: New Phenomenon and Diffusion Model-Augmented Construction}                      

\author[1,2]{Shuyue Fu}[orcid=0009-0001-5111-9779,style=chinese]
\ead{fushuyue@buaa.edu.cn}
\credit{Data curation, Formal analysis, Methodology, Software, Writing - Original draft preparation, Writing - review $\&$ editing}
\affiliation[1]{organization={School of Astronautics, Beihang University},
                addressline={Xueyuan Road No.37}, 
                postcode={100191}, 
                city={Beijing},
                country={People's Republic of China}}

\affiliation[2]{organization={Shen Yuan Honors College, Beihang University},
                addressline={Xueyuan Road No.37}, 
                city={Beijing},
                postcode={100191}, 
                country={People's Republic of China}}

\affiliation[3]{organization={State Key Laboratory of High-Efficiency Reusable Aerospace Transportation Technology}, 
                postcode={102206}, 
                city={Beijing},
                country={People's Republic of China}}

\affiliation[4]{organization={Key Laboratory of Spacecraft Design Optimization $\&$ Dynamic Simulation Technologies, Ministry of Education}, 
                postcode={100191}, 
                city={Beijing},
                country={People's Republic of China}}

\author[1]{Wenxuan Zhang}[style=chinese]
\ead{zwx03511@buaa.edu.cn}
\credit{Methodology, Formal analysis, Writing - review $\&$ editing}

\author[1,4]{Di Wu}[style=chinese]
\ead{wudi2025@buaa.edu.cn}
\credit{Conceptualization, Methodology, Software, Writing - review $\&$ editing}

\author[1,3]{Shengping Gong}[style=chinese]
\cormark[1]
\ead{gongsp@buaa.edu.cn}
\credit{Conceptualization, Funding acquisition, Writing - review $\&$ editing}

\author[1,4]{Peng Shi}[style=chinese]
\ead{shipeng@buaa.edu.cn}
\credit{Methodology, Formal analysis, Writing - review $\&$ editing}

\cortext[cor1]{Corresponding author}

\begin{abstract}
Construction of Earth-Moon transfers is the basis of missions to explore the Moon and cislunar space. The traditional grid search method suffers from a relatively low convergence rate and computational efficiency, mainly focusing on the distribution of transfer characteristic parameters. Moreover, when constructing transfers with different departure/arrival orbital altitudes, the process of grid search and trajectory correction should be repeated with a low convergence rate and computational efficiency. To address these limitations of the traditional grid search method, this paper is devoted to exploring an effective way to augment the grid search method. Bi-impulsive Earth-Moon transfers from a circular Earth parking orbit to a circular Moon target orbit in the Earth-Moon planar circular restricted three-body problem are considered in this paper. Firstly, the transfers are constructed, and the corresponding solution space is explored in terms of construction parameters, including departure phase angle at the Earth parking orbit, initial-to-circular velocity ratio, and time of flight. An interesting phenomenon about the discontinuous behavior of the time-of-flight distribution with respect to departure phase angle is identified. This phenomenon is further used to train a diffusion model, which aims to augment the traditional grid search method and generate high-quality initial guesses for transfers with different departure/arrival orbital altitudes. The construction results and the preliminary applicable range of the proposed method are presented and analyzed. The proposed diffusion model-augmented grid search method improves the convergence rate by 47.34-56.25$\%$ and saves wall-clock time by 39.39-40.52$\%$ with respect to the traditional grid search method, while ensuring comparable transfer characteristics.
\end{abstract}

\begin{keywords}
Planar restricted three-body problem \sep Bi-impulsive Earth-Moon transfer \sep Diffusion model \sep Grid search method
\end{keywords}

\maketitle
\onehalfspacing

\section{Introduction}\label{sec1}
Interest in the Earth-Moon transfers has been renewed due to the proposal of several missions to explore the Moon and cislunar space. For these missions, several transfer scenarios have been proposed, including the Earth-Moon transfers (transfers to circular Moon orbit \cite{topputo2013optimal,qi2017optimal,oshima2019low} and multi-body periodic orbits around the Moon, such as distant retrograde orbit \cite{demeyer2007transfer,yin2023midcourse,wang2025mechanism} and distant prograde orbit \cite{mingotti2010mixed,mingotti2012transfers}) and libration point transfers (transfers between/to periodic orbits of collinear libration points \cite{Parker2014,du2022transfer,OSHIMA2025644} and triangular libration points \cite{ZHANG20152899,CAPDEVILA20181826,YIN2026449}). Among these scenarios, the bi-impulsive scenario from a circular Earth parking orbit to a circular Moon target orbit is a classical scenario for scholars to consider \cite{topputo2013optimal}. To construct such transfers, an appropriate dynamical model should be adopted. Several dynamical models have been adopted, including the patched two-body problem \cite{battin1999introduction}, the Earth-Moon planar three-body problem (PCR3BP) \cite{bib49,guo2019families,jing2026study}, the Sun-Earth/Moon planar bicircular restricted four-body problem (PBCR4BP) \cite{topputo2013optimal,qi2017optimal,oshima2019low}, and other higher-fidelity models \cite{belbruno1993sun,wang2024low,campana2025ephemeris}. Compared to transfers constructed in the patched two-body problem \cite{battin1999introduction}, there is a potential to use the multi-body problems (e.g., the Earth-Moon PCR3BP and the Sun-Earth/Moon PBCR4BP) to obtain more transfer families, including direct and low-energy transfers, to satisfy specific mission requirements. Therefore, to investigate more transfer families and further explore the solution space of bi-impulsive Earth-Moon transfers, one of the multi-body problems, i.e., the Earth-Moon PCR3BP, is adopted as the dynamical model to construct transfer trajectories.

Because there is no closed-form solution in the Earth-Moon PCR3BP, numerical methods are adopted to construct bi-impulsive transfers. Typically, the construction method can be divided into two steps: generation of initial guesses and trajectory correction. Firstly, the initial guesses of transfer trajectories are generated by several methods, including grid search \cite{topputo2013optimal,oshima2019low} and methods based on the dynamical structures \cite{belbruno1993sun,Koon2001,fu2025four,hiraiwa2026design}. Then, initial guesses are corrected to satisfy the specific constraints, i.e., the departure/arrival orbital altitudes and tangential departure/arrival \cite{topputo2013optimal}. When adopting the methods based on the dynamical structures, invariant manifolds \cite{Koon2001}, weak stability boundary structures \cite{belbruno1993sun,fu2026deep}, sequences of lobe dynamics \cite{hiraiwa2026design}, and other dynamical structures using multi-body prior knowledge \cite{oshima2017analysis,qi2017optimal} are usually used. These methods can effectively obtain low-energy transfers, but suffer from a limited number of obtained solutions \cite{dutt2018review}. In contrast, the grid search method performs an exhaustive search of construction parameters and obtains feasible solutions using trajectory correction. This type of method typically yields more solutions and provides a global map to analyze the transfer characteristics \cite{topputo2013optimal,oshima2019low}. Topputo \cite{topputo2013optimal} and Oshima et al. \cite{oshima2019low} presented a global map of impulse and time of flight (TOF) for the bi-impulsive transfers within 100 and 200 Days in the Sun-Earth/Moon PBCR4BP, revealing the rich structure of the solution space. Based on their work, Campana and Topputo \cite{campana2024clustering} used the data mining method to extract the transfer families in the Sun-Earth/Moon PBCR4BP. The aforementioned work \cite{topputo2013optimal,oshima2019low,campana2024clustering} performed a valuable exploration on the solution space and transfer characteristics of bi-impulsive Earth-Moon transfers in the multi-body problems. However, they did not focus on the distribution of construction parameters. Such analysis can further provide useful insight into selecting high-quality initial guesses and aid in the construction of bi-impulsive transfers. Pinelli \cite{pinelli2023neural} pioneered this type of exploration using the neural network method to learn the mapping between construction parameters and generate the initial guesses for trajectory correction accordingly. However, Pinelli \cite{pinelli2023neural} mainly focused on one type of transfer trajectory, i.e., exterior transfer without fly-bys. Therefore, a more systematic exploration of the solution space in terms of construction parameters should be performed. Combining with other limitations of the grid search method, we summarize the limitations of the current work on the grid search method to construct Earth-Moon transfers in the multi-body problems:

\begin{enumerate}
  \item Generally, using the grid search method to generate initial guesses can be considered as a cold-start initialization, which typically yields a relatively low convergence rate.
  \item Also, because the grid search method is a cold-start initialization, the corresponding trajectory correction usually requires a long computational time.
  \item Current work \cite{topputo2013optimal,oshima2019low,campana2024clustering} on the grid search method typically focus on transfer characteristics and transfer families, but pay less attention to the distribution of construction parameters, which can provide further insight into trajectory construction.
  \item Current work \cite{topputo2013optimal,oshima2019low,campana2024clustering} on the grid search method typically focus on transfers with the specific departure/arrival orbital altitudes (i.e., transfers from a 167 km circular Earth parking orbit to a 100 km circular Moon target orbit). When considering transfers with other altitudes, the process of grid search and trajectory correction should be repeated, with a relatively low convergence rate and computational efficiency.
\end{enumerate}

To address the aforementioned limitations, \textbf{the main purpose of this paper is to explore a way to augment the traditional grid search method}. We firstly focus on the solution space of bi-impulsive Earth-Moon transfers with different departure/arrival orbital altitudes, in particular, the distribution of construction parameters (departure phase angle at the Earth parking orbit, initial-to-circular velocity ratio, and time of flight (TOF)). An interesting phenomenon about the discontinuous behavior of the TOF distribution for each departure phase angle is discovered. This phenomenon can also be observed for different settings of orbital altitudes. Furthermore, the distributions of TOF and departure phase angle at the Earth parking orbit exhibit similarity, which can provide further insight into constructing transfers with different departure/arrival orbital altitudes. These phenomena about the distribution of construction parameters are new or less-reported compared to current work. Then, a way to utilize these phenomena to aid in trajectory construction and augment the traditional grid search method is explored. An alternative way is to utilize the generative models, such as variational autoencoder (VAE) \cite{Litteri2026Generation}, flow-based model \cite{dinh2016density} and diffusion model \cite{graebner2025global}. Compared to the neural network method \cite{pinelli2023neural,gao2026rapid}, which are typically tailored to specific transfer families (e.g., direct transfers \cite{gao2026rapid} and exterior transfers without flybys \cite{pinelli2023neural}) and may require the inputs to lie close to the distribution of feasible solutions (when the selected inputs lie outside the distribution of feasible solutions, the prediction may provide an inaccurate initial guess and reduce the convergence rate), generative models learn and approximate the distribution of interest. They select samples from a simple distribution (e.g., a standard Gaussian distribution) and then perform a denoising process to generate data samples. These models have been applied to the initial guess generation of multi-body periodic orbits \cite{litteri2025generative,Litteri2026Generation,clark2026autoencoder}, Earth-Mars ballistic transfers \cite{Presser2024Diffusion}, and low-thrust transfer trajectories \cite{li2023amortized,graebner2025global}. Litteri et al. \cite{litteri2025generative,Litteri2026Generation} and Clark and Scheeres \cite{clark2026autoencoder} used the VAE and autoencoder neural networks to generate the three-body periodic orbits. Similarly, Presser et al. \cite{Presser2024Diffusion} introduced the diffusion model to generate the Earth-Mars transfers. The aforementioned work trained the generative models using the time-discretized trajectory state sequences and did not focus on the distributions of initial states of the considered trajectories/orbits. As a result, when training the generative models, the dynamical information behind the initial state distribution was not utilized effectively. Furthermore, Li et al. \cite{li2023amortized} revealed that clustering structures exist in the solution space of the low-thrust transfer problem in the PCR3BP, and they are affected by the parameter denoted as the maximum allowable thrust. Using these phenomena, Li et al. \cite{li2023amortized} employed a conditional VAE to learn these structures varying with the maximum allowable thrust, generating initial guesses and accelerating the global search for trajectory design. Based on their work, Graebner et al. \cite{graebner2024learning,graebner2025global} further introduced conditional diffusion models to learn the solution-space distributions of low-thrust transfers and the parameter-dependent structures of solutions in costate space. Their work also resulted in the fast generation of initial guesses of the considered problem. Among the generative models, the diffusion model has been widely applied due to fewer constraints on model architecture and enhanced ability to capture complex data distributions \cite{graebner2025global}. Therefore, this paper adopts the diffusion model to aid in generating high-quality initial guesses and augmenting the traditional grid search method. Also, the use of the diffusion model enables the generation of samples with a preset number and is not limited to the original data. This paper combines the dynamical prior knowledge (phenomenon about the distribution of construction parameters) with the diffusion model. Therefore, a diffusion model-augmented grid search method to construct bi-impulsive Earth-Moon transfers with different departure/arrival orbital altitudes is proposed in this paper. Differing from the aforementioned work \cite{li2023amortized,graebner2024learning,graebner2025global}, our newly discovered phenomena motivate us to adopt a shared diffusion model (rather than a conditional diffusion model) to capture the common characteristics of the Earth-Moon transfers with different departure/arrival orbital altitudes. The construction results and the preliminary applicable range of the proposed method are then presented and discussed. Comparison with the results obtained from the traditional grid search method indicates that the proposed method takes advantage of improved convergence rate and computational efficiency for all three settings of departure/arrival orbital altitudes, while ensuring comparable transfer characteristics.

\section{Three-Body Bi-Impulsive Earth-Moon Transfer}\label{sec2}
This section introduces the bi-impulsive Earth-Moon transfers in the Earth-Moon PCR3BP, including the dynamical equations of the Earth-Moon PCR3BP, the concept, and solution space of bi-impulsive transfers in the Earth-Moon PCR3BP. The dynamical analysis of the solution space of bi-impulsive Earth-Moon transfers with different departure/arrival orbital altitudes can provide further insight into the training of the diffusion model.
\subsection{Earth-Moon PCR3BP}\label{subsec2.1}
We consider the construction problem of bi-impulsive Earth-Moon transfers in the Earth-Moon PCR3BP. In this model, the Earth, the Moon, and the spacecraft are assumed to move in the same plane, where the Earth and the Moon move in circular orbit around their barycenter. Selecting the dimensionless units (the length unit (LU) is set as the Earth-Moon distance, the mass unit (MU) is set as the Earth-Moon combined mass, and the time unit (TU) is set as the $1/2\pi$ of the Earth-Moon orbital period around their barycenter) and adopting the Earth-Moon rotating frame, the dynamical equations of the Earth-Moon PCR3BP can be written as
\begin{equation}
\begin{gathered}
  \dot x = u \hfill \\
  \dot y = v \hfill \\
  \dot u = x + 2v - \frac{{\left( {1 - \mu } \right)\left( {x + \mu } \right)}}{{{r_1}^3}} - \frac{{\mu \left( {x + \mu  - 1} \right)}}{{{r_2}^3}} \hfill \\
  \dot v = y - 2u - \frac{{\left( {1 - \mu } \right)y}}{{{r_1}^3}} - \frac{{\mu y}}{{{r_2}^3}} \hfill \\ 
  {r_1} = \sqrt {{{\left( {x + \mu } \right)}^2} + {y^2}} \text{ }\text{ }\text{ }\text{ }{r_2} = \sqrt {{{\left( {x + \mu - 1} \right)}^2} + {y^2}}\hfill \\ 
\end{gathered} \label{eq1}
\end{equation}
where $x$ and $y$ denote the position components with respect to the Earth-Moon barycenter in the Earth-Moon rotating frame, $u$ and $v$ denote the velocity components corresponding to $x$ and $y$ in the Earth-Moon rotating frame, $\mu$ denotes the mass parameter of the Earth-Moon PCR3BP which can be expressed as $\mu=m_\text{Moon}/\left(m_\text{Earth}+m_\text{Moon}\right)$ ($m$ denotes the mass), $r_1$ denotes the distance between the Earth and the spacecraft, and $r_2$ denotes the distance between the Moon and the spacecraft. Numerical integration of Eq. \eqref{eq1} achieved by a variable-step, variable-order (VSVO) Adams-Bashforth-Moulton predictor–corrector method (orders 1 to 13) with absolute and relative tolerances set to $1\times10^{-13}$ is performed to obtain trajectories in the Earth-Moon PCR3BP. The parameters used in this paper can be found in Ref. \cite{topputo2013optimal}, as presented in Table \ref{parameters}. In this paper, bi-impulsive Earth-Moon transfers are constructed in this model. Subsequently, the concept and solution space of bi-impulsive Earth-Moon transfers in this model are introduced.
\begin{table}[h]
\centering
\caption{Parameters of the Earth-Moon PCR3BP}\label{parameters}%
\renewcommand{\arraystretch}{1.5}
\begin{tabular}{@{}llll@{}}
\hline
Symbol & Value  & Units & Meaning\\
\hline
$\mu$    & $1.21506683 \times {10^{ - 2}}$   & --  & Earth-Moon mass parameter  \\
$R_{\text{E}}$    & $6378.145$   & km  & Mean Earth’s radius  \\
$R_{\text{M}}$    & $1737.100$   & km  & Mean Moon’s radius  \\
LU    & $3.84405000 \times {10^5}$   & km  & Length unit  \\
TU    & $3.75676968 \times {10^5}$   & s  & Time unit  \\
\hline
\end{tabular}
\end{table}
\subsection{Concept of Bi-impulsive Earth-Moon Transfers}\label{subsec2.2}
In this paper, we focus on the classical Earth-Moon transfers, particularly the bi-impulsive ones. This type of transfer describes a scenario where the spacecraft is sent to the transfer trajectory from a circular Earth parking orbit with the altitude $h_i$ using an Earth injection impulse $\Delta v_i$, and enters a circular Moon target orbit with the altitude $h_f$ using a Moon insertion impulse $\Delta v_f$. Tangential impulses are performed, and then the states of the transfer trajectories should satisfy the constraints \cite{topputo2013optimal}

\begin{equation}
{\bm{\psi }_i} = \left[ {\begin{array}{*{20}{c}}
  {{{\left( {{x_i} + \mu } \right)}^2} + {y_i}^2 - {{\left( {{R_{\text{Earth}}} + {h_i}} \right)}^2}} \\ 
  {\left( {{x_i} + \mu } \right)\left( {{u_i} - {y_i}} \right) + {y_i}\left( {{v_i} + {x_i} + \mu } \right)} 
\end{array}} \right] = \mathbf{0} \label{eq2}
\end{equation}
\begin{equation}
{\bm{\psi }_f} = \left[ {\begin{array}{*{20}{c}}
  {{{\left( {{x_f} + \mu  - 1} \right)}^2} + {y_f}^2 - {{\left( {{R_{\text{Moon}}} + {h_f}} \right)}^2}} \\ 
  {\left( {{x_f} + \mu  - 1} \right)\left( {{u_f} - {y_f}} \right) + {y_f}\left( {{v_f} + {x_f} + \mu  - 1} \right)} 
\end{array}} \right] = \mathbf{0} \label{eq3}
\end{equation}
The subscript “\textit{i}” denotes quantities corresponding to the epoch when the spacecraft departs from the Earth parking orbit, while the subscript “\textit{f}” denotes quantities corresponding to the epoch when the spacecraft enters the Moon target orbit. To construct transfer trajectories satisfying Eqs. \eqref{eq2}-\eqref{eq3}, we select the construction parameters
\begin{equation}
{\bm{y}} = {\left[ {{\alpha _i},{\text{ }}{\beta _i},{\text{ TOF}}} \right]^{\text{T}}} \label{eq4}
\end{equation}
In these parameters, $\alpha _i$ denotes the departure phase angle of the circular Earth parking orbit in the Earth-Moon rotating frame, $\beta _i$ denotes the initial-to-circular velocity ratio, and TOF denotes the time of flight. With these parameters, the departure states of transfer trajectories can be expressed, following the constraint Eq. \eqref{eq3} as
\begin{equation}
\begin{gathered}
{x_i} = {r_i}\cos {\alpha _i} - \mu \hfill \\
{y_i} = {r_i}\sin {\alpha _i} \hfill \\
{u_i} =  - \left( {{\beta _i}\sqrt {\frac{{1 - \mu }}{{{r_i}}}}  - {r_i}} \right)\sin {\alpha _i}\hfill \\
{v_i} = \left( {{\beta _i}\sqrt {\frac{{1 - \mu }}{{{r_i}}}}  - {r_i}} \right)\cos {\alpha _i}\hfill \\
\label{eq5}
\end{gathered}
\end{equation}
where $r_i=R_\text{Earth}+h_i$. Then, the construction of bi-impulsive Earth-Moon transfers is transformed into searching for the feasible solutions integrated from the departure states shown in Eq. \eqref{eq5} within the time interval $\left[0,\text{ TOF}\right]$ and satisfying the constraint Eq. \eqref{eq3}. This paper explores a diffusion model-augmented method to construct the bi-impulsive Earth-Moon transfers in the Earth-Moon PCR3BP. To provide original data and prior knowledge for the training of the diffusion model, the solution space of the bi-impulsive Earth-Moon transfers is firstly presented and analyzed using the grid search method. 
\subsection{Grid Search Method of Bi-Impulsive Earth-Moon Transfers}\label{subsec2.3}
In this paper, three groups of $h_i$ and $h_f$ are considered: 
\begin{equation} \begin{gathered}
  {\text{Group I:     }}{h_i} = 167{\text{ km,     }}{h_f} = 100{\text{ km}} \hfill \\
  {\text{Group II:    }}{h_i} = 167{\text{ km,     }}{h_f} = 1500{\text{ km}} \hfill \\
  {\text{Group III:   }}{h_i} = 1500{\text{ km,   }}{h_f} = 100{\text{ km}} \hfill \\
\end{gathered} 
\label{eq6}
\end{equation}
To generate initial guesses, we set $\bm{y}$ in Eq. \eqref{eq4} as ${\alpha _i} \in \left[ {0,{\text{ }}2\pi } \right){\text{ rad}}$ with a step-size of $\pi /36{\text{ rad}}$, ${\beta _i} \in \left[ {1.4,{\text{ }}1.414} \right]$ with a step-size of 0.0001, and ${\text{TOF}} \in \left[ {\frac{\pi}{30}\frac{\text{TU}}{86400},{\text{ 8}}\pi \frac{\text{TU}}{86400} } \right]{\text{ Day}}$ with a step-size of $\frac{\pi}{30}\frac{\text{TU}}{86400}{\text{ Day}}$. The ranges of $\alpha_i$ and $\beta_i$ are selected according to Ref. \cite{oshima2019low}, and the selection of TOF does not affect the qualitative findings reported in this paper. With these settings, we obtain 2436480 initial guesses for each group in total. Then, the initial guess is generated from the initial state calculated by Eq. \eqref{eq5}.  

Once the initial guesses are generated, trajectory correction is performed to make the trajectories satisfy the constraints Eqs. \eqref{eq2}-\eqref{eq3}. Trajectory correction can be transformed into a nonlinear programming (NLP) problem. When solving the NLP problem, the sequential quadratic programming (SQP) method is adopted. Since the generated initial guesses satisfy the constraint Eq. \eqref{eq2} rigorously, the constraint of the NLP problem is set as Eq. \eqref{eq3}, and the objective function is set to 0. During the trajectory correction, the ranges of the NLP variables are set as ${\alpha _i} \in \left( { - \infty ,{\text{ }}\infty } \right){\text{ rad}}$, ${\beta _i} \in \left[ { 1.4 ,{\text{ }}1.414 } \right]$, and ${\text{TOF}} \in \left[ {\frac{\pi}{30}\frac{\text{TU}}{86400},{\text{ 8}}\pi \frac{\text{TU}}{86400} } \right]{\text{ Day}}$ (during trajectory correction, the unit of TOF is the dimensionless unit TU). To ensure the accuracy and efficiency of the SQP method, the parameters used in the trajectory correction are determined through trial and error. The tolerance of the step-size during the iteration is set to $1\times 10^{-12}$, the tolerance of the constraint violation is set to $1 \times 10^{-8}$, and the first-order optimality tolerance is set to $1 \times 10^{-6}$. Moreover, the maximum number of both iterations and function evaluations are set to 1000. During solving the NLP problem, the Earth/Moon collision trajectories are excluded. When the correction is finished, the solutions satisfying $\left\| {{\bm{\psi }_f}} \right\| < 1 \times {10^{ - 8}}$ are selected and recorded. The values of $\alpha_i$ are then transformed into those within $\left[0,\text{ }2\pi\right]{\text{ rad}}$.
\subsection{Solution Space of Bi-impulsive Earth-Moon Transfers}\label{subsec2.4}
\begin{figure}
\centering
\includegraphics[width=0.99\textwidth]{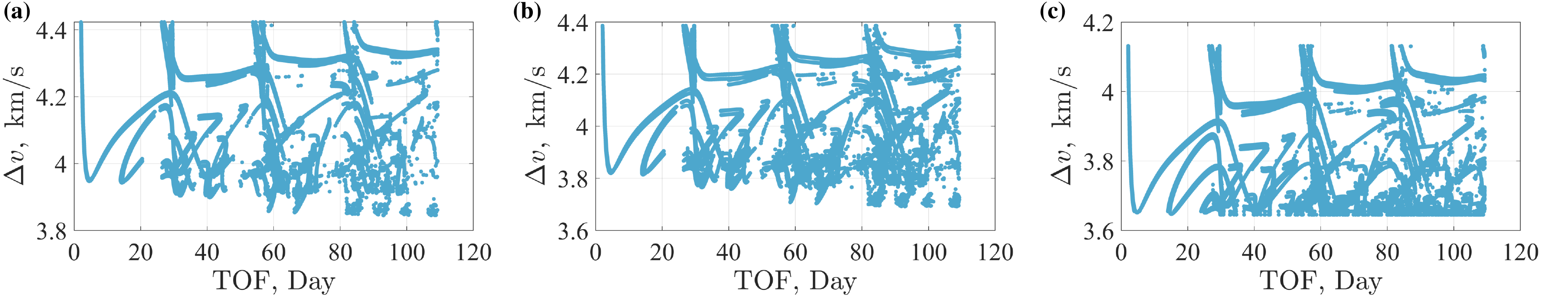}
\caption{The $\left( {{\text{TOF}},{\text{ }}\Delta v} \right)$ map of obtained bi-impulsive Earth-Moon transfers in the Earth-Moon PCR3BP. (a) Group I; (b) Group II; (c) Group III.}
\label{fig_map}
\end{figure}

Using the aforementioned method, we totally obtain 1018001 feasible solutions for Group I, 1132795 solutions for Group II, and 1212226 solutions for Group III. The corresponding solution space is presented in the terms of the $\left( {{\text{TOF}},{\text{ }}\Delta v} \right)$ map shown in Fig. \ref{fig_map}. This type of map has been widely adopted \cite{topputo2013optimal,oshima2019low} to describe the solution space of the transfers in the multi-body problem. However, although the map provides the full information about the transfer characteristics, it does not provide more information about the distribution of the construction parameters. In this paper, we would like to present the distribution of the construction parameters and report an interesting, novel, or less-reported phenomenon, i.e., the discontinuous behavior of TOF distribution.

\begin{figure}
\centering
\includegraphics[width=0.99\textwidth]{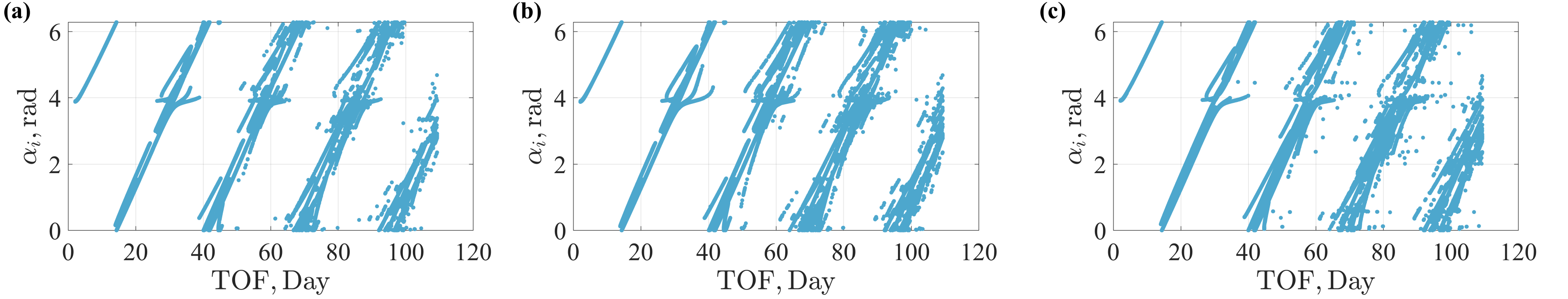}
\caption{The $\left( {{\text{TOF}},{\text{ }}\alpha _i} \right)$ map of obtained bi-impulsive Earth-Moon transfers in the Earth-Moon PCR3BP. (a) Group I; (b) Group II; (c) Group III.}
\label{fig_map_TOF_alphai}
\end{figure}

\begin{figure}
\centering
\includegraphics[width=0.99\textwidth]{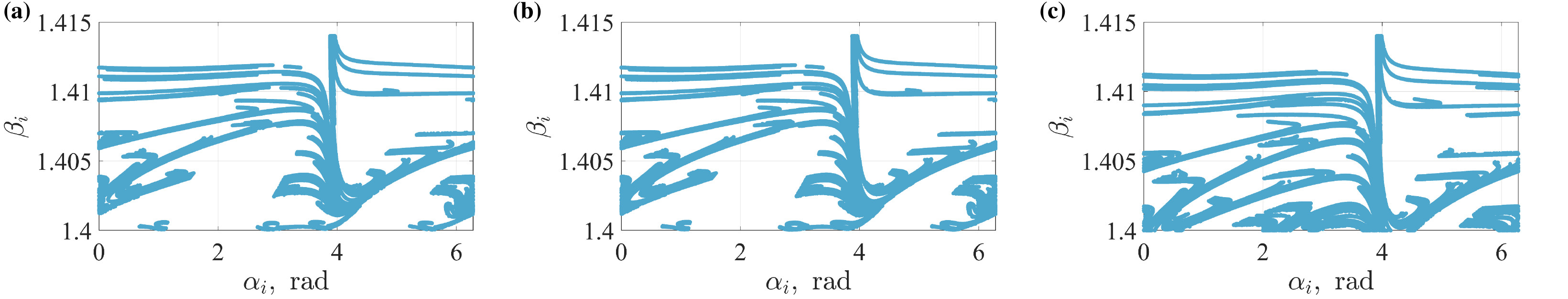}
\caption{The $\left( \alpha_i,{\text{ }}\beta _i \right)$ map of obtained bi-impulsive Earth-Moon transfers in the Earth-Moon PCR3BP. (a) Group I; (b) Group II; (c) Group III.}
\label{fig_map_alphai_betai}
\end{figure}

For the obtained solutions, the $\left( {{\text{TOF}},{\text{ }}{\alpha _i}} \right)$ map and the $\left( {{\alpha _i},{\text{ }}{\beta _i}} \right)$ map for three groups are presented in Figs. \ref{fig_map_TOF_alphai}-\ref{fig_map_alphai_betai}. The construction parameters for these four groups exhibit similar distributions. Differing from the $\left( {{\alpha _i},{\text{ }}{\beta _i}} \right)$ map, the $\left( {{\text{TOF}},{\text{ }}{\alpha _i}} \right)$ map reveals a rather pronounced phenomenon:
\begin{enumerate}
  \item For each $\alpha_i$, the distribution of TOF reveals a discontinuous behavior. The map can be approximately divided into 5 branches, and the TOF of each branch (central locations) differs by approximately one month.
  \item For each branch, $\alpha_i$ exhibits a banded distribution with respect to TOF, and the slope of $\alpha_i$ with respect to TOF is approximately equal.
\end{enumerate}

\begin{figure}
\centering
\includegraphics[width=0.66\textwidth]{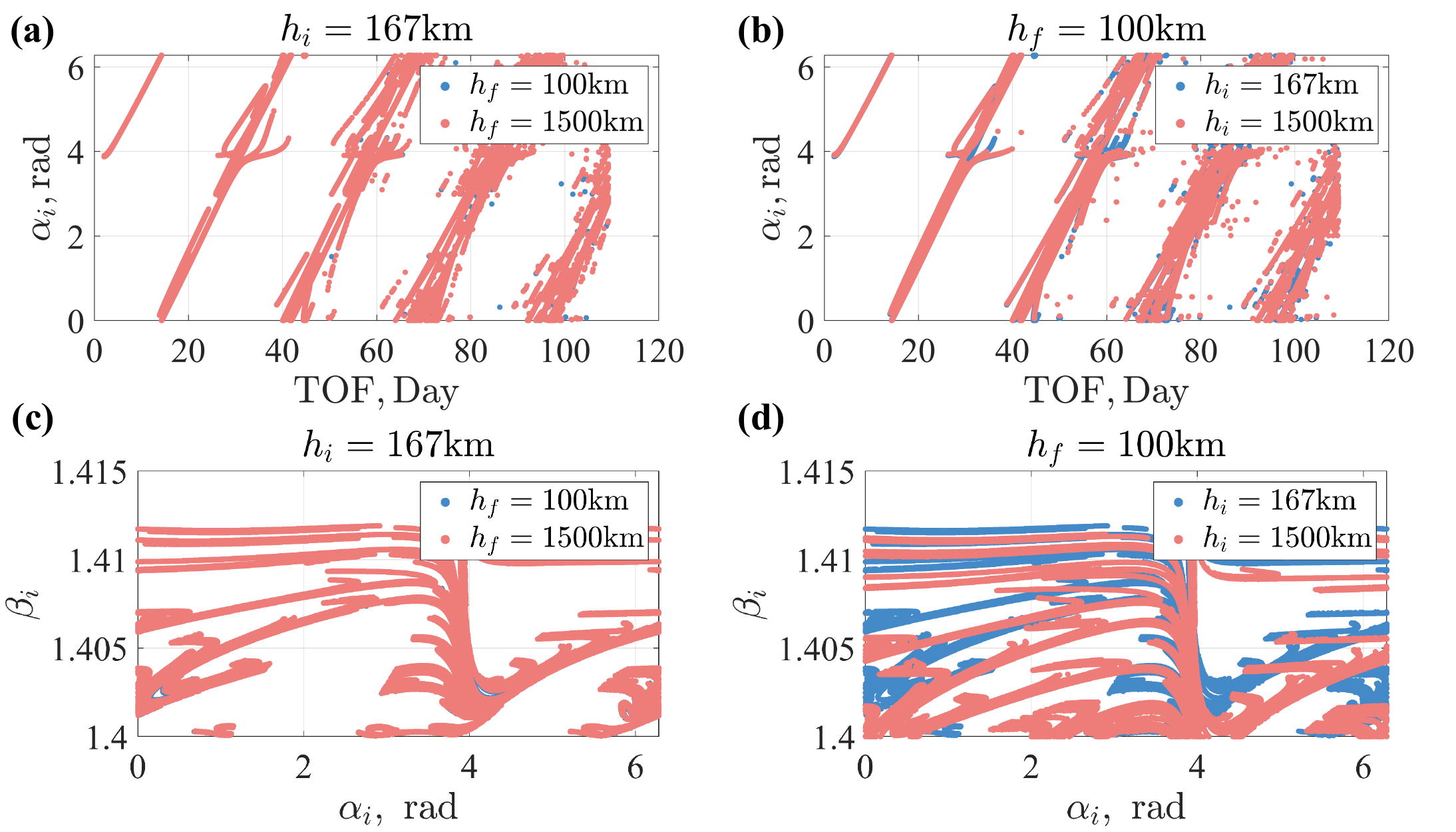}
\caption{Solution space for different groups. (a) $\left( {{\text{TOF}},{\text{ }}\alpha _i} \right)$ map for Groups I and II; (b) $\left( {{\text{TOF}},{\text{ }}\alpha _i} \right)$ map for Groups I and III; (c) $\left( {{\alpha _i},{\text{ }}{\beta _i}} \right)$ map for Groups I and II; (d) $\left( {{\alpha _i},{\text{ }}{\beta _i}} \right)$ map for Groups I and III.}
\label{different_group}
\end{figure}

The aforementioned phenomenon exists in all of three groups, implying its robustness with respect to the orbital altitude settings. Also, when adopting a smaller step-size of construction parameters, a similar phenomenon has been observed for Group I \cite{fu2025discontinuous}, implying its robustness with respect to the grid search settings. This phenomenon possibly corresponds to the orbital period of the Earth/Moon around their barycenter, and cannot be revealed by the $\left( {{\text{TOF}},{\text{ }}\Delta v} \right)$ map, as the distribution shown in Fig. \ref{fig_map} reveals a continuity of TOF. Here, we provide a preliminary interpretation for this phenomenon. The detailed interpretation can be found in Appendix A. \textbf{It can further provide useful insight into selecting initial guesses to construct transfers, i.e., when performing the grid search, there might be no need to select the values of $\left( {{\text{TOF}},{\text{ }}{\alpha _i}} \right)$ pair located in the blank region in Fig. \ref{fig_map_TOF_alphai} (a)-(c).} Also, examining the solution space between different groups (see Fig. \ref{different_group}), it can be observed that the $\left( {{\text{TOF}},{\text{ }}\alpha _i} \right)$ maps are similar between solutions with different values of $h_i$ of $h_f$, while the $\left( {{\alpha _i},{\text{ }}{\beta _i}} \right)$ maps can be different more pronouncedly. To quantitatively evaluate the similarity between the $\left( {{\text{TOF}},{\text{ }}\alpha _i} \right)$ maps with different departure/arrival orbital altitudes, the Sliced Wasserstein-1 distance (SW1) metric calculated using POT 0.9.7 in a Python 3.12.12 environment is adopted. The departure phase angle $\alpha_i$ is transformed into $\cos \alpha_i$ and $\sin \alpha_i$. TOF is normalized using the maximum value of $8\pi \frac{\text{TU}}{86400}$ Day and the minimum value of $\frac{\pi}{30} \frac{\text{TU}}{86400}$ Day, and $\cos \alpha_i$ and $\sin \alpha_i$ are normalized using the maximum value of 1 and the minimum value of $-1$. A lower SW1 metric implies a higher degree of similarity between the two distributions. In this paper, we adopt 2048 deterministic projection directions (the random seed is set to 42) to calculate SW1 between three groups of $\left( {{\text{TOF}},{\text{ }}\alpha _i} \right)$ maps. The results are presented as follows: the SW1 metric between Groups I and II is 0.01236, and the SW1 metric between Groups I and III is 0.03694. From these values, it can be observed that the $\left( {{\text{TOF}},{\text{ }}{\alpha _i}} \right)$ maps between Groups I and II exhibit a higher degree of similarity than that between Groups I and III. However, we still think the $\left( {{\text{TOF}},{\text{ }}{\alpha _i}} \right)$ maps between Groups I and III exhibit similarity in the sample distribution and branch width. For comparison, we add $\pi$ to $\alpha _i$ in Group I and then transform it to $\left[0,\text{ }2\pi\right]$. We further calculate the SW1 metric between the $\left( {{\text{TOF}},{\text{ }}\alpha _i} \right)$ map in Group I and this additional “$\left( {{\text{TOF}},{\text{ }}{\alpha _i}} \right)$ maps”, and obtain a value of 0.11499. This value further indicates that the $\left( {{\text{TOF}},{\text{ }}{\alpha _i}} \right)$ maps between Groups I and III also exhibit similarity. This phenomenon motivates us to use the $\left( {{\text{TOF}},{\text{ }}\alpha _i} \right)$ maps to train a diffusion model, as once the model is trained, it can be applied to the construction of bi-impulsive Earth-Moon transfers with different departure/arrival orbital altitudes. Subsequently, a diffusion model-augmented construction method is proposed. The construction results will further support the claim of this similarity.

\section{Diffusion Model-Augmented Construction Method of Bi-Impulsive Earth-Moon Transfers}\label{sec3}
In this section, a diffusion model-augmented method to construct bi-impulsive Earth-Moon transfers mentioned in Section \ref{sec2} is proposed. Firstly, the dataset for training is prepared based on the aforementioned dynamical analysis. Then, the training and evaluation of the diffusion model are presented. Finally, a construction method of bi-impulsive Earth-Moon transfers combining the samples generated by the diffusion model with the grid search is proposed, enabling the generation of high-quality initial guesses and improved convergence rate of trajectory correction over the traditional grid search method described in Section \ref{sec2}. The diffusion model is trained in a Python 3.12.12 environment on NVIDIA GeForce RTX 3050 (8 GB VRAM). We set the random seed to 42 to ensure reproducibility.
\subsection{Preparation of the Dataset}\label{subsec3.1}
\begin{figure}
\centering
\includegraphics[width=0.64\textwidth]{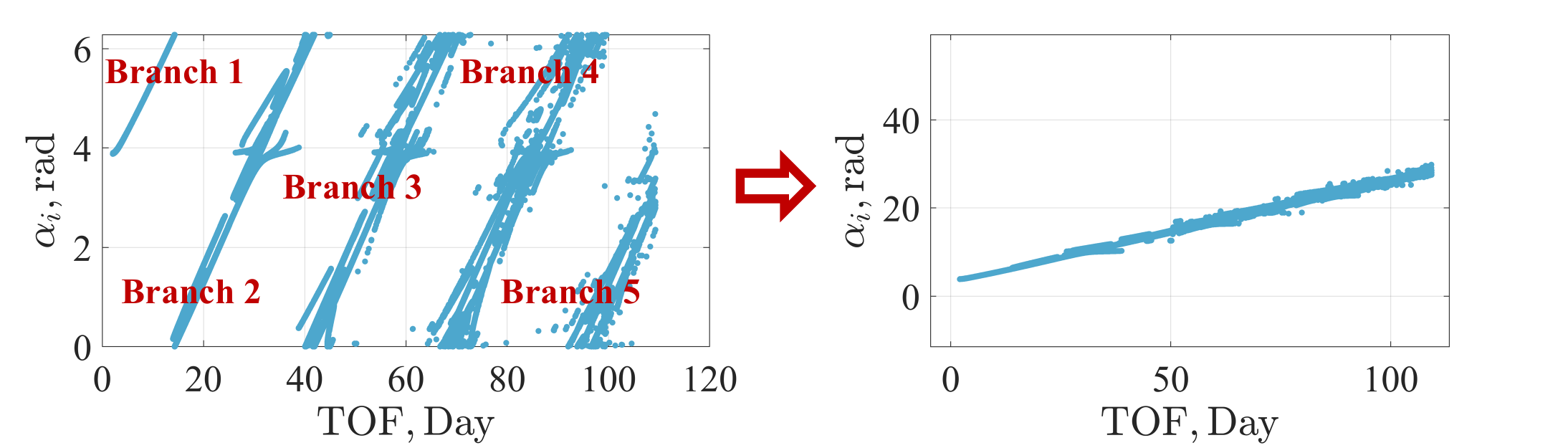}
\caption{The combined $\left( {{\text{TOF}},{\text{ }}\alpha _i} \right)$ map for Group I.}
\label{combine_shiyi}
\end{figure}

Based on the dynamical analysis of the solution space mentioned in Section \ref{sec2}, we adopt the $\left( {{\text{TOF}},{\text{ }}\alpha _i} \right)$ map to train the diffusion model. In particular, only the $\left( {{\text{TOF}},{\text{ }}\alpha _i} \right)$ map with $h_i=167\text{ km}$ and $h_f=100\text{ km}$ (i.e., Group I) is selected due to the similar distributions between different groups shown in Fig. \ref{different_group} (a) and (b). As shown in Fig. \ref{combine_shiyi}, we combine the 5 branches into one branch to simplify the training, using the mapping

\begin{equation} \begin{gathered}
  {\alpha _{i{\text{Branch }}1}} = {\alpha _{i{\text{Branch }}1}}{\text{, TO}}{{\text{F}}_{{\text{Branch }}1}} = {\text{TO}}{{\text{F}}_{{\text{Branch }}1}} \hfill \\
  {\alpha _{i{\text{Branch 2}}}} = {\alpha _{i{\text{Branch 2}}}}{\text{ + 2}}\pi {\text{, TO}}{{\text{F}}_{{\text{Branch 2}}}} = {\text{TO}}{{\text{F}}_{{\text{Branch 2}}}} \hfill \\
  {\alpha _{i{\text{Branch 3}}}} = {\alpha _{i{\text{Branch 3}}}}{\text{ + 4}}\pi {\text{, TO}}{{\text{F}}_{{\text{Branch 3}}}} = {\text{TO}}{{\text{F}}_{{\text{Branch 3}}}} \hfill \\
  {\alpha _{i{\text{Branch 4}}}} = {\alpha _{i{\text{Branch 4}}}}{\text{ + 6}}\pi {\text{, TO}}{{\text{F}}_{{\text{Branch 4}}}} = {\text{TO}}{{\text{F}}_{{\text{Branch 4}}}} \hfill \\
  {\alpha _{i{\text{Branch 5}}}} = {\alpha _{i{\text{Branch 5}}}}{\text{ + 8}}\pi {\text{, TO}}{{\text{F}}_{{\text{Branch 5}}}} = {\text{TO}}{{\text{F}}_{{\text{Branch 5}}}} \hfill \\ 
\end{gathered} 
\label{eq7}
\end{equation}

Once the dataset (i.e., the transformed $\left( {{\text{TOF}},{\text{ }}\alpha _i} \right)$ map shown in Fig. \ref{combine_shiyi}) is obtained, we use the \textit{ShuffleSpilt} function from sklearn to divide the dataset into the training dataset (80$\%$) and the validation dataset (20$\%$). Then, we denote the transformed $\left( {{\text{TOF}},{\text{ }}\alpha _i} \right)$ points as $\bm{x}_0$, and perform a normalization of these data:
\begin{equation} \tilde{\bm{x}}_0 = (\bm{x}_0-\bm{\mu}_\text{Training}) \oslash\left(\bm{\sigma}_\text{Training}+1\times 10^{-6}\right)
\label{eq8}
\end{equation}
where $\bm{\mu}_\text{Training}$ and $\bm{\sigma}_\text{Training}$ denote the mean and standard deviation of the \textbf{training dataset}, and the subscript “Training” denotes quantities corresponding to the training dataset. When generating samples, the inverse normalization is performed. Subsequently, the training of the diffusion model is presented.
\subsection{Diffusion Model}\label{subsec3.2}
\begin{figure}
\centering
\includegraphics[width=0.69\textwidth]{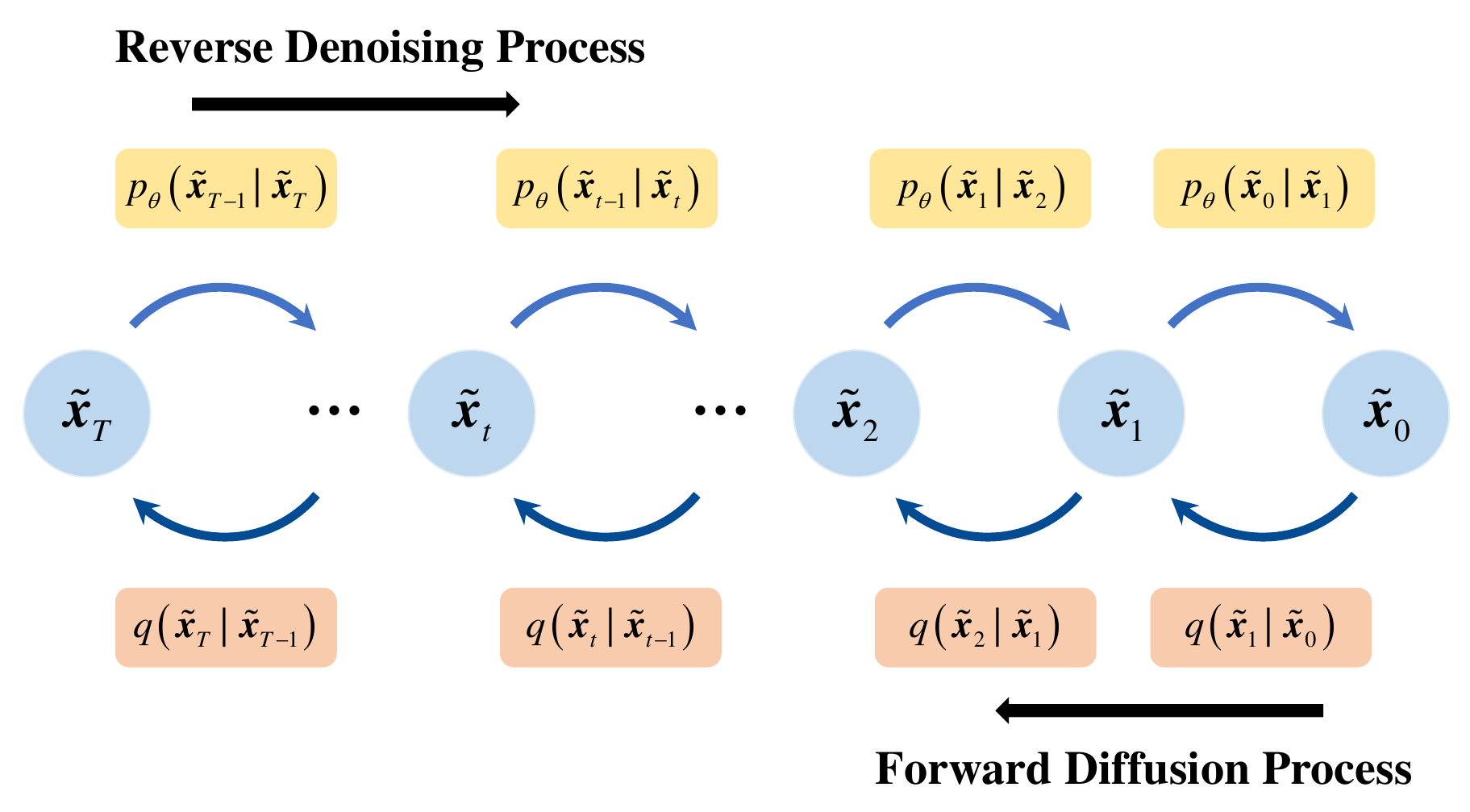}
\caption{Schematic of the forward diffusion process and reverse denoising process.}
\label{diffusion_model}
\end{figure}

When using the diffusion model to generate samples, the samples from a standard Gaussian distribution are firstly selected. Subsequently, these samples are mapped into samples corresponding to a distribution of interest, i.e., the expected samples \cite{graebner2025global} (generated $\left( {{\text{TOF}},{\text{ }}\alpha _i} \right)$ samples in this paper). This process is denoted as the reverse denoising process. Moreover, to train a diffusion model, the forward diffusion process is required. Therefore, the diffusion model can typically be divided into two processes, i.e., the forward diffusion process and reverse denoising process (shown in Fig. \ref{diffusion_model}). Then, the details about these processes are introduced, and the training of the diffusion model is presented.
\subsubsection{Forward Diffusion Process and Reverse Denoising Process}\label{subsubsec3.2.1}
The forward diffusion process describes a scenario where the Gaussian noise is gradually added into the samples $\tilde {\bm{x}}_0$, and finally makes the samples approximated as an approximately standard Gaussian noise, i.e., $q\left( {{{\tilde {\bm{x}}}_T}} \right) \approx \mathcal{N}\left( {{{\tilde {\bm{x}}}_T};{\text{ }}\bm{0},{\text{ }}\bm{I}} \right)$, where $T$ denotes the total diffusion step and $\bm{I}$ is the identity matrix. To achieve this target, $q\left( {{{\tilde {\bm{x}}}_t}|{{\tilde {\bm{x}}}_{t - 1}}} \right)$ in the forward diffusion process is presented as (in Fig. \ref{diffusion_model}, $t$ denotes an intermediate diffusion step) 
\begin{equation} q\left( {{{\tilde {\bm{x}}}_t}|{{\tilde {\bm{x}}}_{t - 1}}} \right) = \mathcal{N}\left( {{{\tilde {\bm{x}}}_t};{\text{ }}\sqrt {{a_t}} {{\tilde {\bm{x}}}_{t - 1}},{\text{ }}{b_t}{\bm{I}}} \right), \text{ }a_t=1-b_t
\label{eq9}
\end{equation}
where $b_t$ denotes the noise variance at diffusion step $t$, and $\{b_t\}_{t=1}^T$ denotes the variance schedule. In this paper, the linear variance schedule
\begin{equation} {b_t} = {b_{\min }} + \frac{{t - 1}}{{T - 1}}\left( {{b_{\max }} - {b_{\min }}} \right)
\label{eq10}
\end{equation}
is adopted \cite{NEURIPS2020_4c5bcfec}, where $b_{\min}=1\times 10^{-4}$, $b_{\max}=0.02$, and $T=1000$. With Eq. \eqref{eq9}, ${\tilde {\bm{x}}_t}$ can be constructed explicitly from ${\tilde {\bm{x}}_0}$:
\begin{equation} {\tilde {\bm{x}}_t} = \sqrt {{{\bar a}_t}} {\tilde {\bm{x}}_0} + \sqrt {1 - {{\bar a}_t}} \bm{\epsilon} ,{\text{ }}\bm{\epsilon}  \sim \mathcal{N}\left(\bm{\epsilon};\text{ } {\bm{0},{\text{ }}\bm{I}} \right)
\label{eq11}
\end{equation}
where:
\begin{equation} {\bar a_t} = \prod\limits_{s = 1}^t {{a_s}}
\label{eq12}
\end{equation}
In the reverse denoising process, samples from the standard Gaussian distribution are selected, and the noise ${\bm{\epsilon} _{\bm{\theta}} }\left( {{{\tilde {\bm{x}}}_t},{\text{ }}t} \right) $ is removed iteratively. During this process, a neural network with parameters $\bm{\theta}$ is used to predict the noise ${\bm{\epsilon} _{\bm{\theta}} }\left( {{{\tilde {\bm{x}}}_t},{\text{ }}t} \right) $. In Fig. \ref{diffusion_model}, ${p_{\bm{\theta}} }\left( {{{\tilde {\bm{x}}}_{t - 1}}|{{\tilde {\bm{x}}}_t}} \right)$ is written as
\begin{equation} {p_{\bm{\theta}} }\left( {{{\tilde {\bm{x}}}_{t - 1}}|{{\tilde {\bm{x}}}_t}} \right) = \mathcal{N}\left( {{{\tilde {\bm{x}}}_{t - 1}};{\text{ }}{{\bm{\mu}} _{\bm{\theta}} }\left( {{{\tilde {\bm{x}}}_t},{\text{ }}t} \right),{\text{ }}{{\overset{\lower0.5em\hbox{$\smash{\scriptscriptstyle\frown}$}}{b} }_t}\bm{I}} \right),\text{ }{{\overset{\lower0.5em\hbox{$\smash{\scriptscriptstyle\frown}$}}{b} }_t}=b_t \frac{1-\bar a_{t-1}}{1-\bar a_t}
\label{eq13}
\end{equation}
where
\begin{equation} {\bm{\mu} _{\bm{\theta}} }\left( {{{\tilde {\bm{x}}}_t},{\text{ }}t} \right) = \frac{1}{{\sqrt {{a_t}} }}\left( {{{\tilde {\bm{x}}}_t} - \frac{{{b_t}}}{{\sqrt {1 - {{\bar a}_t}} }}{\bm{\epsilon} _{\bm{\theta}} }\left( {{{\tilde {\bm{x}}}_t},{\text{ }}t} \right)} \right)
\label{eq14}
\end{equation}
Therefore, the predicted ${\tilde {\bm{x}}_{t - 1}}$ can be obtained as
\begin{equation} {\tilde {\bm{x}}_{t - 1}} = \left\{ {\begin{array}{*{20}{c}}
  {{\bm{\mu} _{\bm{\theta}} }\left( {{{\tilde {\bm{x}}}_t},{\text{ }}t} \right) + \sqrt {{{\overset{\lower0.5em\hbox{$\smash{\scriptscriptstyle\frown}$}}{b} }_t}} \bm{z},{\text{ }}\bm{z} \sim \mathcal{N}\left(\bm{z};\text{ } {\bm{0},{\text{ }}\bm{I}} \right),{\text{ }}t > 1} \\ 
  {{\bm{\mu} _{\bm{\theta}} }\left( {{{\tilde {\bm{x}}}_t},{\text{ }}t} \right),{\text{ }}t = 1} 
\end{array}} \right.
\label{eq15}
\end{equation}
Based on the aforementioned framework, samples corresponding to a distribution of interest can be generated, which can consequently provide high-quality initial guesses for bi-impulsive Earth-Moon transfer construction. Subsequently, the training and evaluation of the diffusion model in this paper are presented.
\subsubsection{Training and Evaluation of the Diffusion Model}\label{subsubsec3.2.3}
\begin{figure}
\centering
\includegraphics[width=0.95\textwidth]{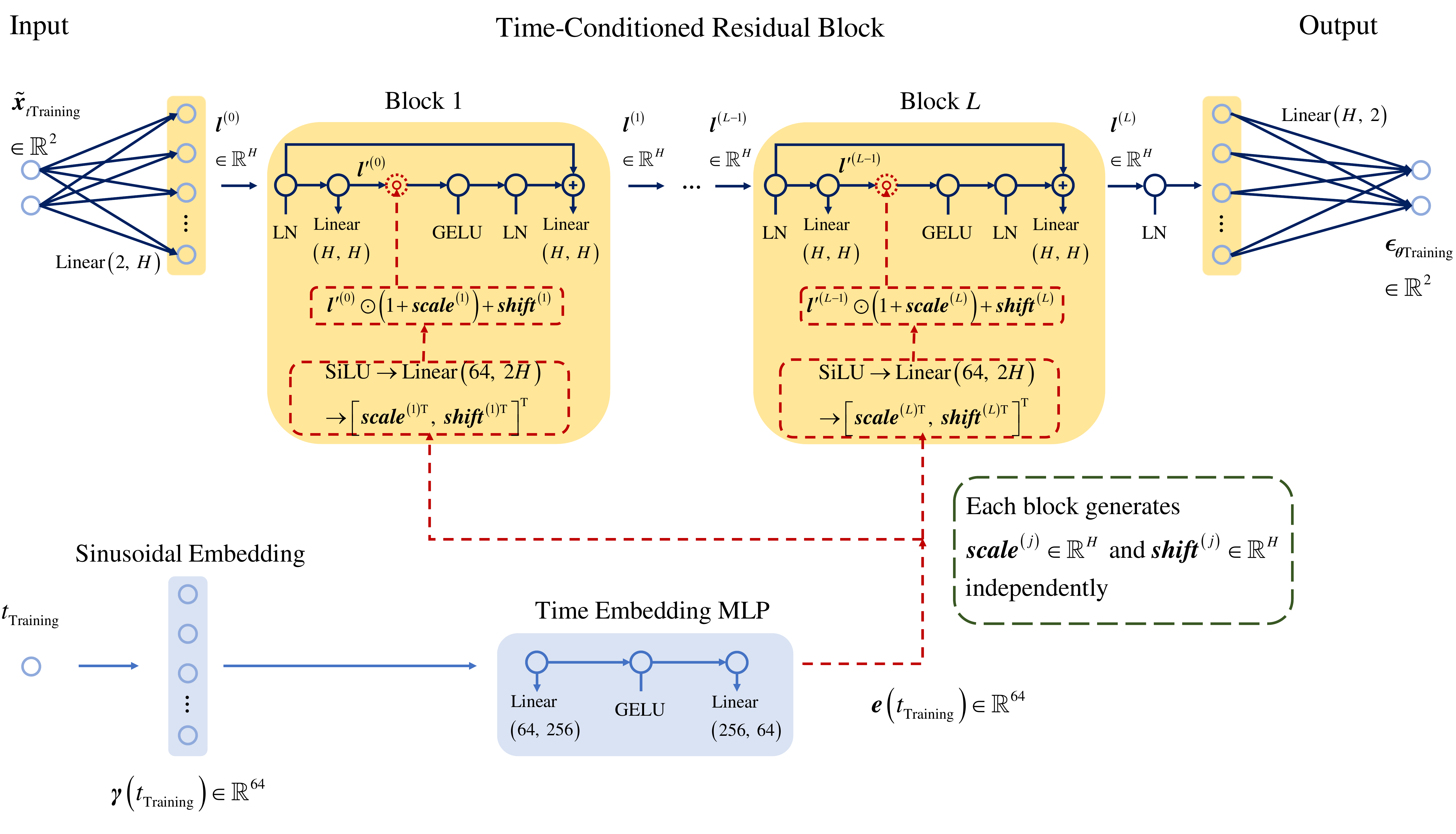}
\caption{Schematic of the structure of the neural network.}
\label{NN}
\end{figure}

Using the training dataset and validation dataset mentioned in Section \ref{subsec3.1}, the training and evaluation of the diffusion model are presented. The diffusion model is trained using PyTorch 2.9.0 (CUDA 13.0). In this paper, a time-conditioned residual multilayer perceptron (MLP) is adopted as the neural network to predict the noise, whose structure is shown in Fig. \ref{NN}. The input of the neural network is set as $\tilde {\bm{x}}_{t\text{Training}}$ and diffusion step $t_\text{Training}$. In particular, the diffusion step $t$ is transformed into a 64-dimensional time embedding using a sinusoidal embedding \cite{NIPS2017_3f5ee243} followed by a time embedding ($64\to 256 \to 64$) MLP. The resulting time embedding is further projected in each residual block to generate FiLM scale and shift parameters, as shown in Fig. \ref{NN}. The sinusoidal embedding can be written as
\begin{equation} \begin{gathered}
  \bm{\gamma }\left( {{t_{{\text{Training}}}}} \right) \in {\mathbb{R}^{64}} \hfill \\
  \bm{\gamma }\left( {{t_{{\text{Training}}}}} \right) = \left[ {\sin \left( {{t_{{\text{Training}}}}{\omega _0}} \right),\text{ }...,\text{ }\sin \left( {{t_{{\text{Training}}}}{\omega _{31}}} \right)},\text{ } {\cos \left( {{t_{{\text{Training}}}}{\omega _0}} \right),\text{ }...,\text{ }\cos \left( {{t_{{\text{Training}}}}{\omega _{31}}} \right)}\right] \hfill \\
  {\omega _i} = \frac{1}{{{{10000}^{\frac{i}{{64/2 - 1}}}}}},\text{ }i=0,\text{ }1,\text{ }...,\text{ }31 \hfill \\ 
\end{gathered} 
\label{eq20}
\end{equation}

The output of the neural network is set as $\bm{\epsilon}_{\bm{\theta}}\left(\tilde {\bm{x}}_{t\text{Training}},\text{ }t_\text{Training}\right)$. In Fig. \ref{NN}, the term “LN" denotes layer normalization. When training the neural network to predict the noise, the batch size is fixed as 1024. For each mini-batch, $\tilde {\bm{x}}_{0\text{Training}}$ are selected from the training dataset using a shuffled \textit{DataLoader}. For each $\tilde {\bm{x}}_{0\text{Training}}$, $t_\text{Training}$ and $\bm{\epsilon}_\text{Training}$ are sampled:
\begin{equation} \begin{gathered}
  t_\text{Training} \sim \mathcal{U}\left\{ {1,{\text{ }}...,\text{ }T} \right\} \hfill \\
  \bm{\epsilon}_\text{Training}  \sim \mathcal{N}\left( {\bm{\epsilon}_\text{Training} ;{\text{ }}\bm{0},{\text{ }}\bm{I}} \right) \hfill \\ 
\end{gathered}
\label{eq16}
\end{equation}
Then, $\tilde {\bm{x}}_{t\text{Training}}$ can be obtained by Eq. \eqref{eq11}. Then, the neural network predicts the predicted noise $\bm{\epsilon}_{\bm{\theta}}\left(\tilde {\bm{x}}_{t\text{Training}},\text{ }t_\text{Training}\right)$. We adopt the mean squared error (MSE) between $\bm{\epsilon}_\text{Training}$ and $\bm{\epsilon}_{\bm{\theta}}\left(\tilde {\bm{x}}_{t\text{Training}},\text{ }t_\text{Training}\right)$ as the loss function \cite{NEURIPS2020_4c5bcfec} evaluated by \textit{mse\_loss} function
\begin{equation} \text{Loss}_{\text{Training}}\left(\bm{\theta}\right)=\mathbb{E}_{\tilde{\bm{x}}_{0\text{Training}},t_\text{Training},\bm{\epsilon}_\text{Training}}\left[\text{MSE}\left(\bm{\epsilon}_\text{Training},\text{ }\bm{\epsilon}_{\bm{\theta}}\left(\tilde{\bm{x}}_{t\text{Training}},\text{ }t_\text{Training}\right)\right)\right]
\label{eq17}
\end{equation}
Based on the aforementioned discussion, the training framework of the used neural network can be summarized as follows:
\begin{enumerate}
  \item \textit{Structure of the neural network:} A time-conditioned residual MLP is adopted as the neural network to predict the noise, and the schematic of its structure is shown in Fig. \ref{NN}. For this structure, the number of time-conditioned residual blocks ($L$) and hidden dimension ($H$) are treated as the hyperparameters (see Table \ref{tab2}).
\item \textit{Training process:} The \textit{Adam} optimizer is employed, and gradients are clipped with a maximum norm of 1.0. The loss function is set as Eq. \eqref{eq17}, and the learning rate ($lr$) is treated as the hyperparameter (see Table \ref{tab2}). The batch size is fixed as 1024. For each mini-batch iteration, a mini-batch estimate of the loss function shown in Eq. \eqref{eq17} is calculated and the parameters $\bm{\theta}$ of the neural network are updated. The training epoch is set to 500. When recording the loss function corresponding to the training dataset for each epoch, the weighted average with respect to each mini-batch is adopted.
\end{enumerate}

\begin{table}[h]
\centering
\caption{Settings of Hyperparameters}\label{tab2}%
\begin{tabular}{@{}ll@{}}
\hline
Hyperparameter & Setting   \\
\hline
$L$    & 4, 8, 16 \\
$H$    & 64, 128, 256 \\
$lr$ & $1\times 10^{-3}$, $1\times 10^{-4}$, $1\times 10^{-5}$ \\
\hline
\end{tabular}
\end{table}

When evaluating the performance of the neural network and selecting the optimal hyperparameter combination, the validation dataset is employed. For each validation mini-batch, $\tilde {\bm{x}}_{0\text{Validation}}$ are selected using a \textit{DataLoader} which is not shuffled. The Monte Carlo validation is performed, and $t_\text{Validation}$ and $\bm{\epsilon}_\text{Validation}$ are sampled 8 times (i.e., $t_\text{Validation}^{\left(k\right)}$ and $\bm{\epsilon}_\text{Validation}^{\left(k\right)}$, $k=1,\text{ }...,\text{ }8$) to calculate the loss function corresponding to the validation dataset:
\begin{equation} 
\text{Loss}_{\text{Validation}}\left(\bm{\theta}\right)=\frac{1}{8}\sum\limits_{k = 1}^8\mathbb{E}_{\tilde{\bm{x}}_{0\text{Validation}},t_\text{Validation}^{\left(k\right)},\bm{\epsilon}_\text{Validation}^{\left(k\right)}}\left[\text{MSE}\left(\bm{\epsilon}_\text{Validation}^{\left(k\right)},\text{ }\bm{\epsilon}_{\bm{\theta}}\left(\tilde{\bm{x}}_{t\text{Validation}}^{\left(k\right)},\text{ }t_\text{Validation}^{\left(k\right)}\right)\right)\right]
\label{eq21}
\end{equation}
During the training, we record the minimum $\text{Loss}_{\text{Validation}}(\bm{\theta})$ and the corresponding epoch for each hyperparameter combination. We select the optimal hyperparameter combination with the minimum $\text{Loss}_{\text{Validation}}(\bm{\theta})$ among these 27 combinations and the corresponding best epoch ($E_\text{best}$). Then, we combine the training dataset and validation dataset, and perform a retraining. The epoch of retraining is set as $\max\left(1,\text{ }\left\lfloor E_{\mathrm{best}}\frac{N_{\text{Training}}+N_{\text{Validation}}}{N_{\text{Training}}}\right\rfloor\right)$ (where $N$ denotes the sample number). Once the retraining is completed, the $\left( {{\text{TOF}},{\text{ }}\alpha _i} \right)$ samples can be generated.
\subsubsection{Generation of New Samples}\label{subsubsec3.2.4}
When using the trained diffusion model to generate new samples, we select 20736 points from the standard Gaussian distribution
\begin{equation} 
{\tilde {\bm{x}}_{T{\text{Generated}}}} \sim \mathcal{N}\left( {{{\tilde {\bm{x}}}_{T{\text{Generated}}}};{\text{ }}\bm{0},{\text{ }}\bm{I}} \right)
\label{eq22}
\end{equation}
Then, we use the trained neural network to predict the noise and perform a reverse denoising process shown in Eq. \eqref{eq15} to obtain $\tilde {\bm{x}}_{0{\text{Generated}}}$. Finally, the inverse normalization is performed:
\begin{equation} 
{{\bm{x}}_{0{\text{Generated}}}} = {\tilde {\bm{x}}_{0{\text{Generated}}}} \odot {\bm{\sigma} _{{\text{Training}}}} + {\bm{\mu} _{{\text{Training}}}}
\label{eq23}
\end{equation}
Subsequently, with ${{\bm{x}}_{0{\text{Generated}}}}$, we propose a diffusion model-augmented grid search method to construct bi-impulsive Earth-Moon transfers with different departure/arrival orbital altitudes.

\subsection{Diffusion Model-Augmented Grid Search Method}\label{subsec3.3}
To improve the convergence rate and computational efficiency of the traditional grid search method, we use the aforementioned diffusion model to augment this type of method. The main steps of the proposed diffusion model-augmented grid search method are summarized as
\begin{enumerate}
  \item Transform $\alpha_i$ in ${{\bm{x}}_{0{\text{Generated}}}}$ into $\left[0,\text{ }2\pi\right]\text{ rad}$, and select ${\text{TOF}} \in \left[ {\frac{\pi}{30}\frac{\text{TU}}{86400},{\text{ 8}}\pi \frac{\text{TU}}{86400} } \right]{\text{ Day}}$ in ${{\bm{x}}_{0{\text{Generated}}}}$ and the corresponding $\alpha_i$.
  \item Combine the transformed/selected $\left( {{\text{TOF}},{\text{ }}\alpha _i} \right)$ from ${{\bm{x}}_{0{\text{Generated}}}}$ with ${\beta _i} \in \left[ {1.4,{\text{ }}1.414} \right]$ with a step-size of 0.0001 to generate initial guesses. All three groups of altitude settings shown in Eq. \eqref{eq6} are considered.
  \item Perform the trajectory correction following the settings mentioned in Section \ref{subsec2.3}.
\end{enumerate}
Subsequently, the construction results of the proposed method are presented and discussed in Section \ref{sec4}. The advantages over the traditional grid search method are also revealed according to the comparison between these two methods.

\section{Results and Discussion}\label{sec4}
This section presents the results of this paper. Firstly, the loss time history for the aforementioned hyperparameter combinations is presented, which determines the optimal hyperparameter combination with the minimum validation loss. Then, the new samples are generated from the retrained diffusion model, and the diffusion model-augmented grid search method is performed to construct the solution space of bi-impulsive Earth-Moon transfers with different departure/arrival orbital altitudes. Finally, the comparison between the diffusion model-augmented grid search method and the traditional grid search method is performed to reveal the advantages of the proposed method. We perform the trajectory correction on two shared-memory multi-core CPUs (96 cores), and use a loop-level parallelism in the outer loop. Two Intel Xeon Platinum 8488C processors (2.4 GHz base frequency and 256 GB RAM) in a Linux environment are used in the simulations.
\subsection{Selection of the Optimal Hyperparameter Combination}\label{subsec4.1}
\begin{figure}
\centering
\includegraphics[width=0.75\textwidth]{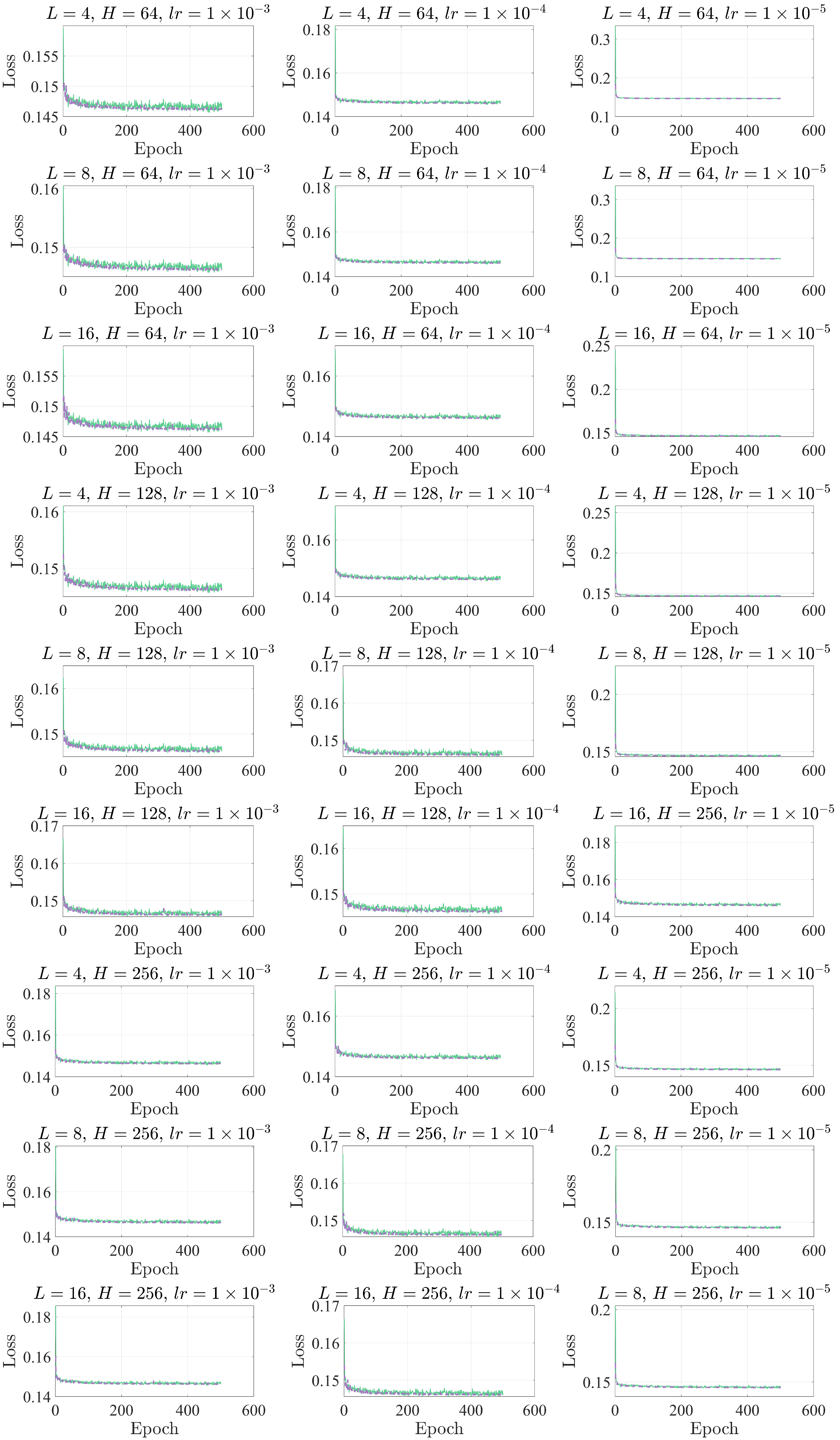}
\caption{The loss time history for the aforementioned hyperparameter combinations.}
\label{fig_loss}
\end{figure}

The loss time history for the aforementioned hyperparameter combinations is presented in Fig. \ref{fig_loss}, where the solid curves denote the training loss while the dashed curves denote the validation loss. From Fig. \ref{fig_loss}, it can be observed that the model does not exhibit obvious overfitting during training. Then, we select the optimal hyperparameter combination with the minimum validation loss and the best epoch.

\begin{figure}
\centering
\includegraphics[width=0.33\textwidth]{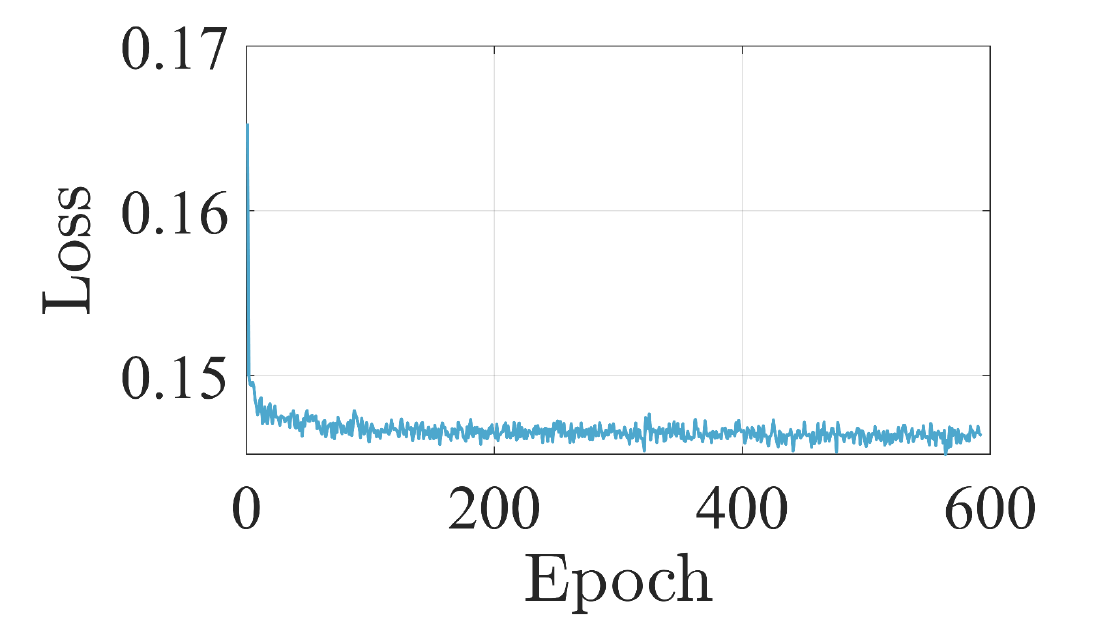}
\caption{The loss time history of retraining.}
\label{fig_retrain}
\end{figure}

\begin{table}[h]
\centering
\caption{Minimum Validation Loss and the Corresponding Best Epoch for the Aforementioned Hyperparameter Combinations}\label{tab3}%
\begin{tabular}{@{}lllll@{}}
\hline
$L$ & $H$ & $lr$ & $\text{Loss}_\text{Validation}$ & $E_\text{best}$  \\
\hline
4 & 64 & $1\times 10^{-3}$ & 0.146100 & 389 \\
4 & 64 & $1\times 10^{-4}$ & 0.146007 & 475 \\
4 & 64 & $1\times 10^{-5}$ & 0.146426 & 481 \\
8 & 64 & $1\times 10^{-3}$ & 0.146057 & 449 \\
8 & 64 & $1\times 10^{-4}$ & 0.146006 & 482 \\
8 & 64 & $1\times 10^{-5}$ & 0.146191 & 496 \\
16 & 64 & $1\times 10^{-3}$ & 0.146089 & 424 \\
16 & 64 & $1\times 10^{-4}$ & 0.146005 & 482 \\
16 & 64 & $1\times 10^{-5}$ & 0.146090 & 499 \\
4 & 128 & $1\times 10^{-3}$ & 0.146103 & 389 \\
4 & 128 & $1\times 10^{-4}$ & 0.145998 & 424 \\
4 & 128 & $1\times 10^{-5}$ & 0.146254 & 481 \\
8 & 128 & $1\times 10^{-3}$ & 0.146093 & 496 \\
8 & 128 & $1\times 10^{-4}$ & 0.146013 & 475 \\
8 & 128 & $1\times 10^{-5}$ & 0.146124 & 499 \\
16 & 128 & $1\times 10^{-3}$ & 0.146099 & 396 \\
8 & 128 & $1\times 10^{-4}$ & 0.146021 & 482 \\
8 & 128 & $1\times 10^{-5}$ & 0.146016 & 500 \\
4 & 256 & $1\times 10^{-3}$ & 0.146127 & 485 \\
4 & 256 & $1\times 10^{-4}$ & 0.145993 & 475 \\
4 & 256 & $1\times 10^{-5}$ & 0.146125 & 500 \\
8 & 256 & $1\times 10^{-3}$ & 0.146105 & 389 \\
8 & 256 & $1\times 10^{-4}$ & 0.146007 & 422 \\
8 & 256 & $1\times 10^{-5}$ & 0.146001 & 482 \\
16 & 256 & $1\times 10^{-3}$ & 0.146098 & 389 \\
16 & 256 & $1\times 10^{-4}$ & 0.146011 & 470 \\
16 & 256 & $1\times 10^{-5}$ & 0.146009 & 475 \\
\hline
\end{tabular}
\end{table}
Table \ref{tab3} presents the minimum validation loss and the corresponding best epoch for each hyperparameter combination. We select the optimal hyperparameter combination with the minimum validation loss in Table \ref{tab3}:
\begin{equation} 
L=4,\text{ }H=256, \text{ }lr=1\times 10^{-4}
\label{eq24}
\end{equation}
Then, with this optimal hyperparameter combination, we combine the training dataset and validation dataset, and retrain the diffusion model following the method mentioned in Section \ref{subsubsec3.2.3}. The loss time history of retraining is presented in Fig. \ref{fig_retrain}.

\begin{figure}
\centering
\includegraphics[width=0.33\textwidth]{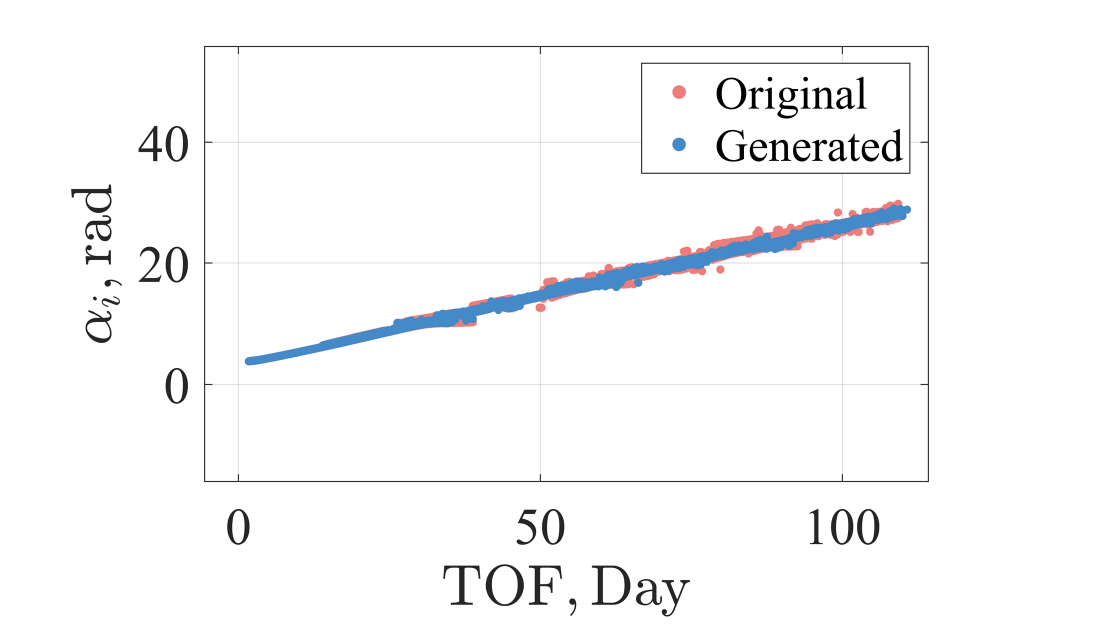}
\caption{The generated samples.}
\label{generated}
\end{figure}

Once the retraining is completed, the new samples are generated, as shown in Fig. \ref{generated}. 
Then, with the generated samples, the proposed diffusion model-augmented grid search method is applied to the construction of bi-impulsive Earth-Moon transfers. Notably, the application of the generated samples is not only limited to constructing Earth-Moon transfers with $h_i=167\text{ km}$ and $h_f=100\text{ km}$, but also with Groups II and III because of the similar $\left( {{\text{TOF}},{\text{ }}\alpha _i} \right)$ distributions shown in Fig. \ref{different_group} (a) and (b). Subsequently, the construction results are presented and discussed.
\subsection{Construction Results of Diffusion Model-Augmented Grid Search Method}\label{subsec4.2}
\begin{figure}
\centering
\includegraphics[width=0.99\textwidth]{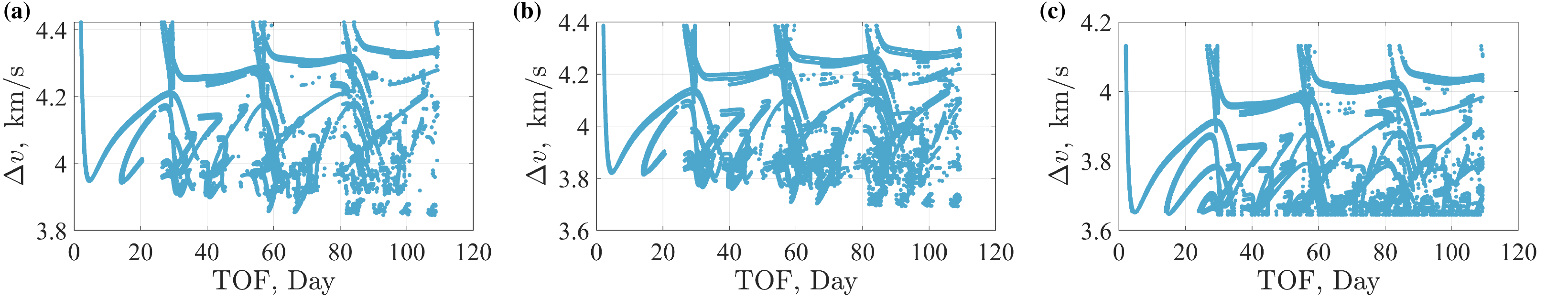}
\caption{The $\left( {{\text{TOF}},{\text{ }}\Delta v} \right)$ map of bi-impulsive Earth-Moon transfers in the Earth-Moon PCR3BP obtained from the diffusion model-augmented grid search method. (a) Group I; (b) Group II; (c) Group III.}
\label{fig_DV_diffusion}
\end{figure}

Following Step 1 in Section \ref{subsec3.3}, we obtain 2920815 initial guesses for Groups I-III, respectively. Figure \ref{fig_DV_diffusion} presents the $\left( {{\text{TOF}},{\text{ }}\Delta v} \right)$ map of bi-impulsive Earth-Moon transfers in the Earth-Moon PCR3BP obtained from the diffusion model-augmented grid search method. Using the diffusion model-augmented grid search method, we totally obtain 1906844 feasible solutions for Group I, 2040593 solutions for Group II, and 2141073 solutions for Group III. These results preliminarily verify the effectiveness of the proposed diffusion model-augmented grid search method.

\begin{figure}
\centering
\includegraphics[width=0.5\textwidth]{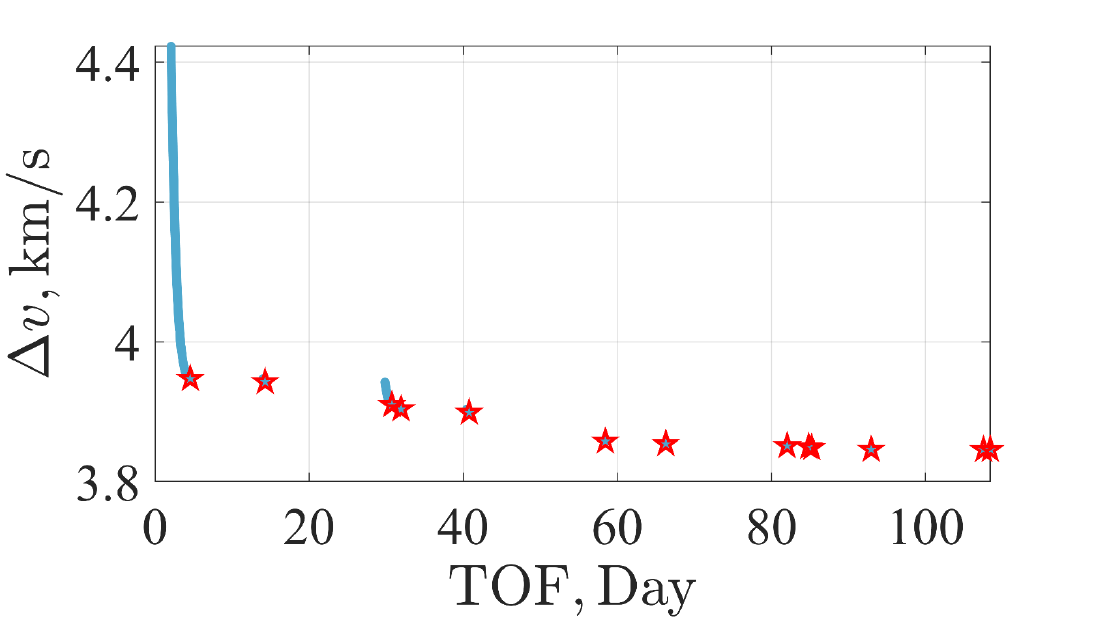}
\caption{Pareto front of the $\left( {{\text{TOF}},{\text{ }}\Delta v} \right)$ map (Group I).}
\label{pareto_group_I}
\end{figure}

\begin{figure}
\centering
\includegraphics[width=0.5\textwidth]{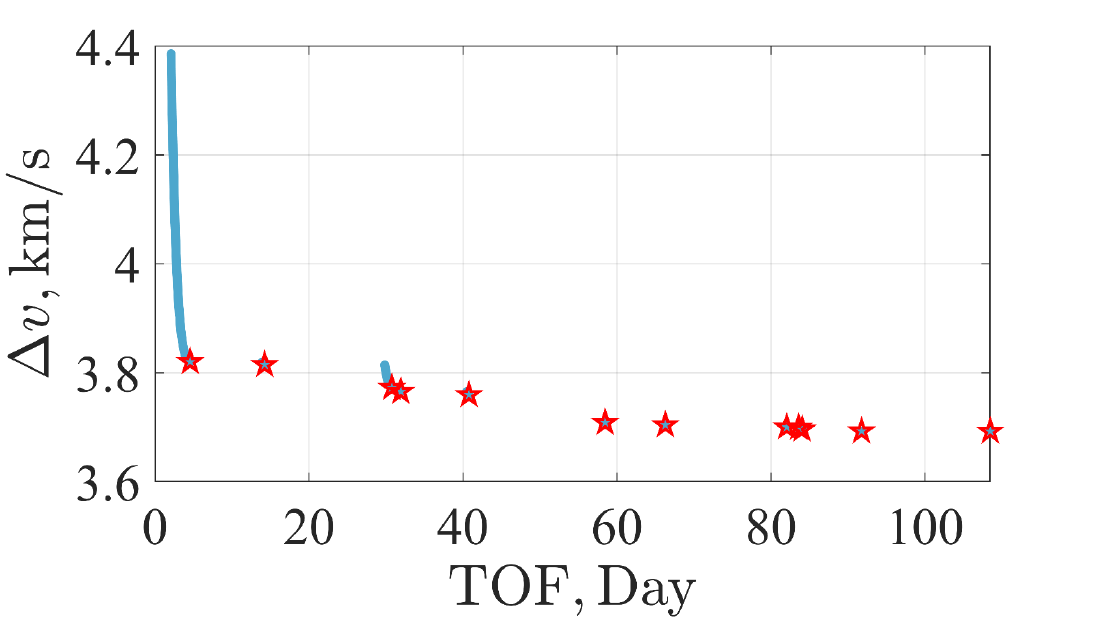}
\caption{Pareto front of the $\left( {{\text{TOF}},{\text{ }}\Delta v} \right)$ map (Group II).}
\label{pareto_group_II}
\end{figure}

\begin{figure}
\centering
\includegraphics[width=0.5\textwidth]{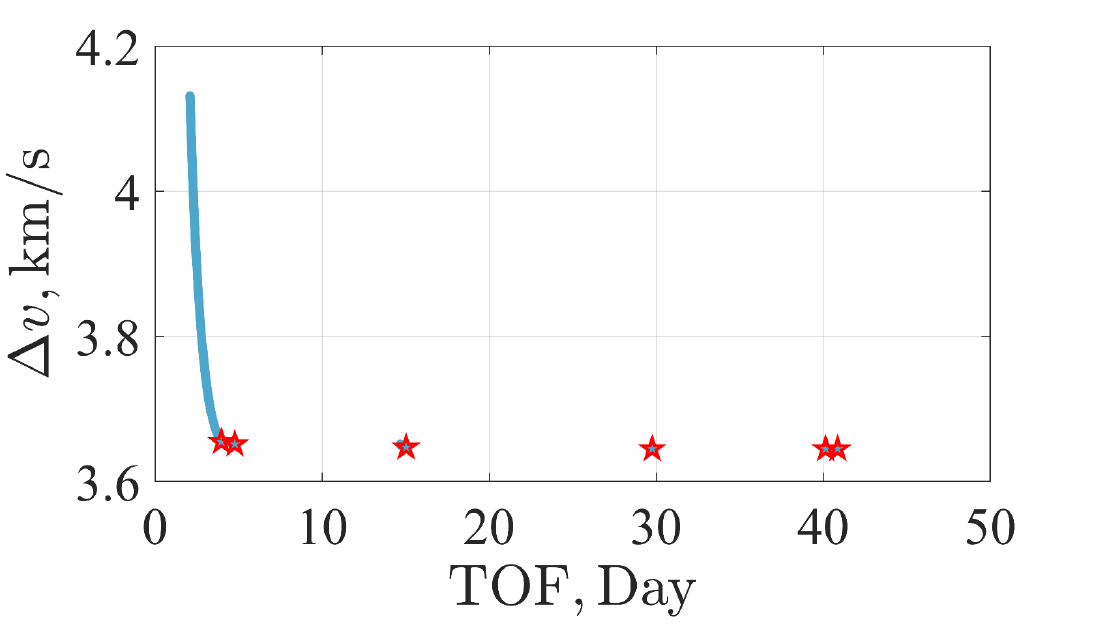}
\caption{Pareto front of the $\left( {{\text{TOF}},{\text{ }}\Delta v} \right)$ map (Group III).}
\label{pareto_group_III}
\end{figure}

\begin{figure}
\centering
\includegraphics[width=0.78\textwidth]{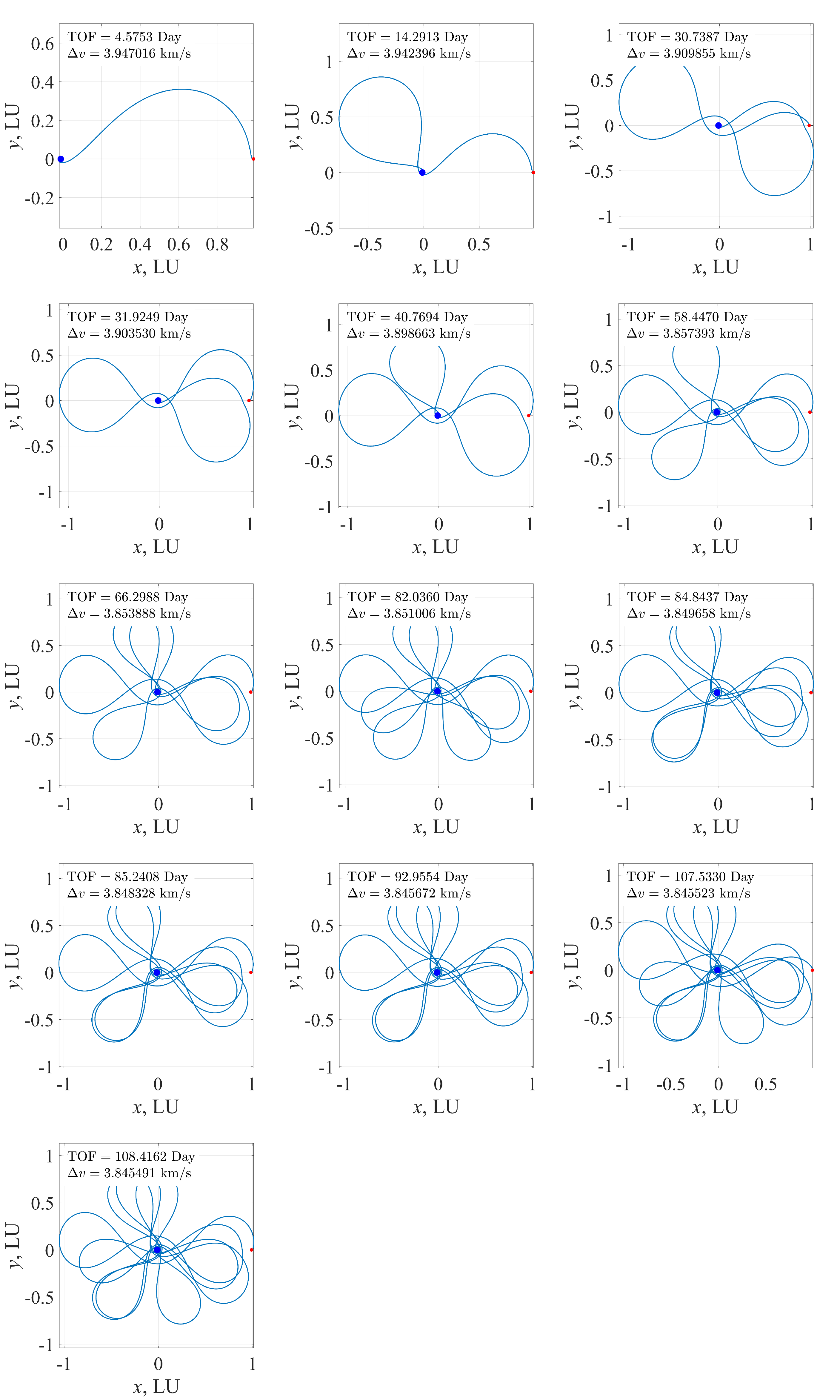}
\caption{Selected trajectory sample (Group I).}
\label{sample_group_I}
\end{figure}

\begin{figure}
\centering
\includegraphics[width=0.78\textwidth]{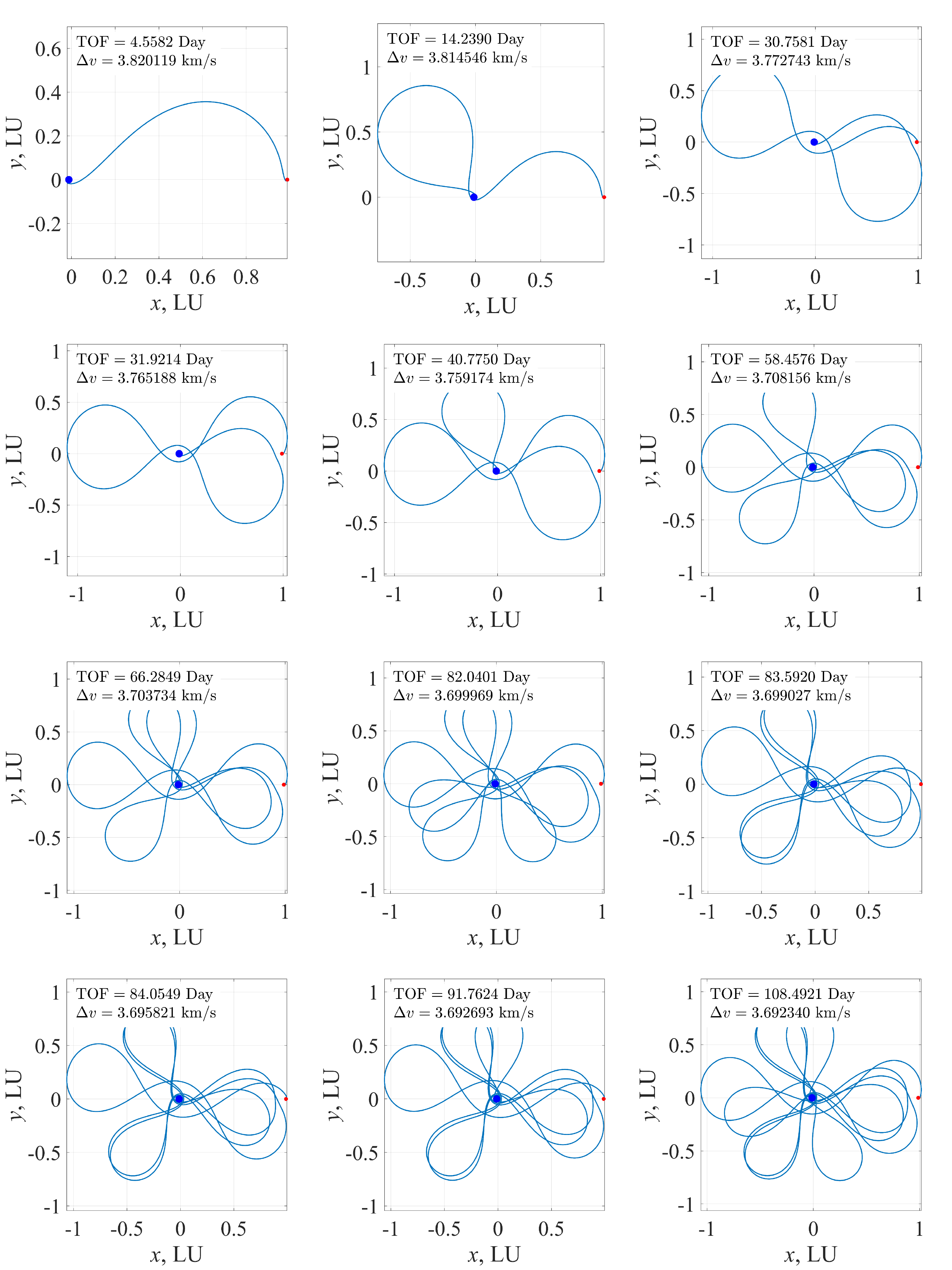}
\caption{Selected trajectory sample (Group II).}
\label{sample_group_II}
\end{figure}

\begin{figure}
\centering
\includegraphics[width=0.78\textwidth]{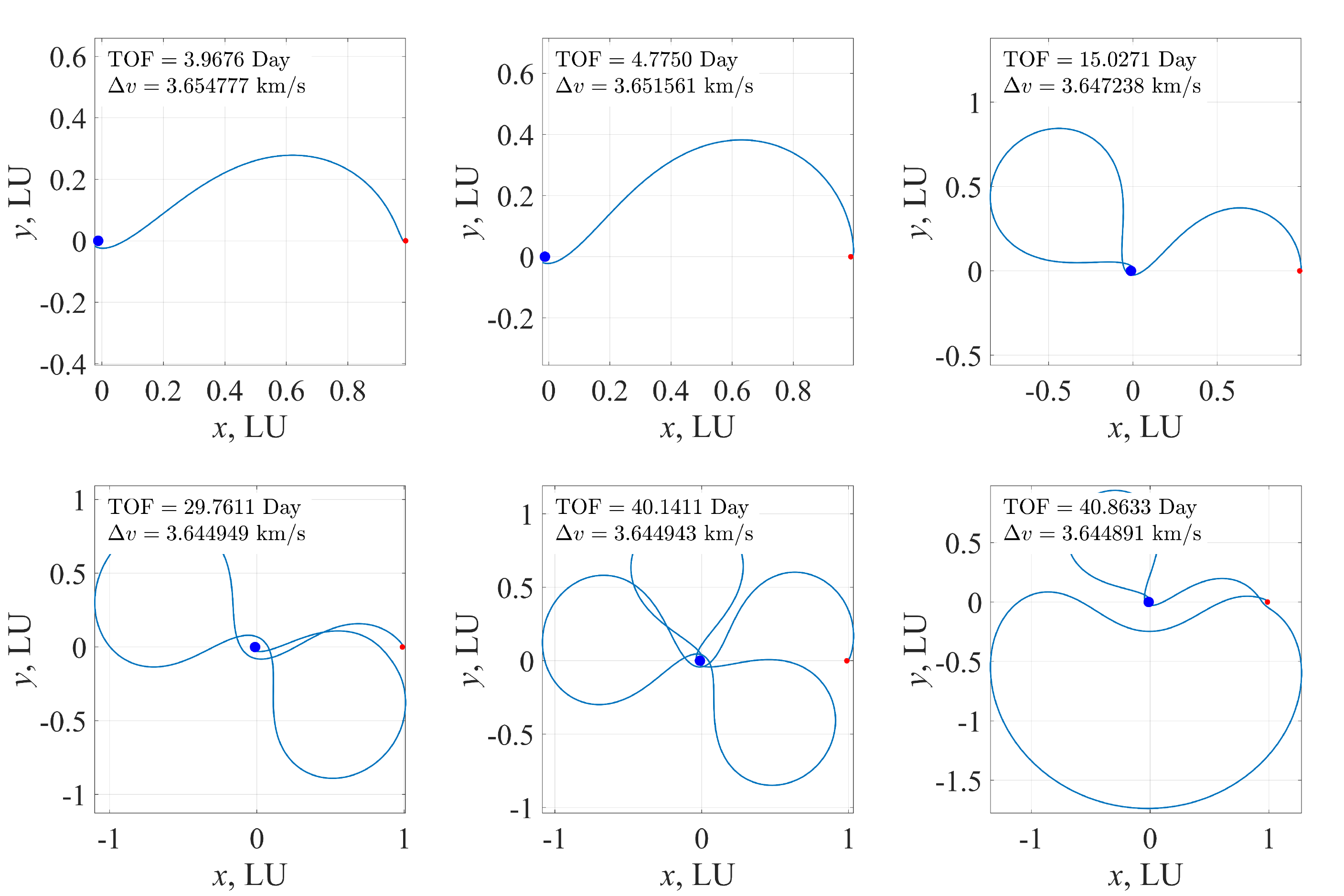}
\caption{Selected trajectory sample (Group III).}
\label{sample_group_III}
\end{figure}

We extract the Pareto front of the $\left( {{\text{TOF}},{\text{ }}\Delta v} \right)$ map shown in Fig. \ref{fig_DV_diffusion}, considering a shorter TOF and lower $\Delta v$, and select typical trajectory samples on the Pareto front. The Pareto fronts for Groups I-III are presented in Figs. \ref{pareto_group_I}-\ref{pareto_group_III}. The red stars denote the corresponding $\left( {{\text{TOF}},{\text{ }}\Delta v} \right)$ of the selected trajectory samples. Figures \ref{sample_group_I}-\ref{sample_group_III} present the typical trajectory samples selected on the Pareto fronts. 

Subsequently, the comparison between results obtained from the proposed diffusion model-augmented grid search method and the traditional grid search method (mentioned in Section \ref{subsec2.3}) is performed to reveal the advantages of the developed method.

\subsection{Comparison with the Traditional Grid Search Method}\label{subsec4.3}
\begin{figure}
\centering
\includegraphics[width=0.75\textwidth]{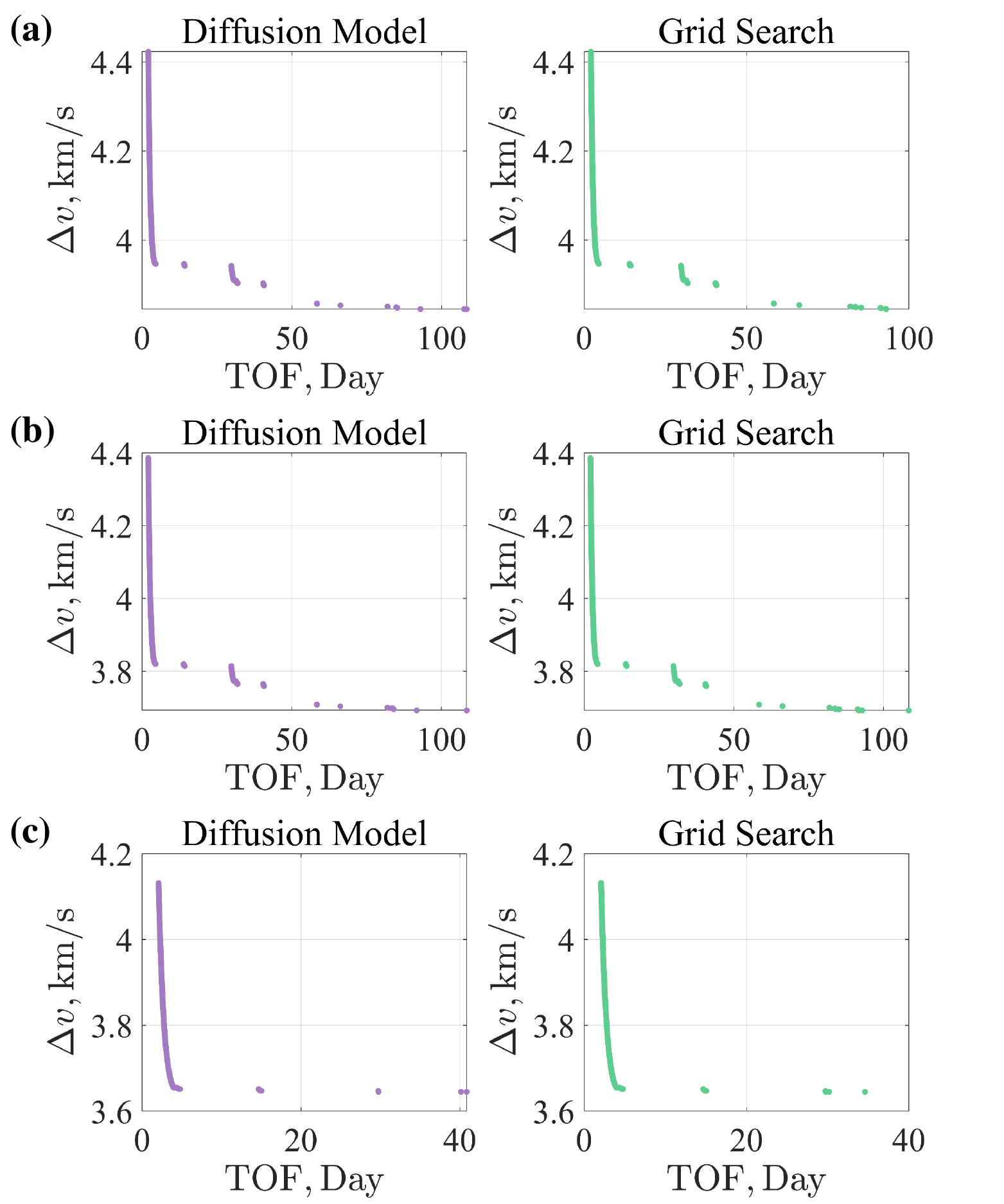}
\caption{Pareto fronts obtained from the two methods. (a) Group I; (b) Group II; (c) Group III.}
\label{comparison_pareto}
\end{figure}

The comparison with the traditional grid search method can be divided into three aspects: transfer characteristics, convergence rate, and computational efficiency. Firstly, a comparison in terms of transfer characteristics is performed. Figure \ref{comparison_pareto} presents the extracted Pareto fronts obtained from these two methods. From Fig. \ref{comparison_pareto}, it can be observed that the Pareto fronts exhibit similarity between the two methods, indicating the similar transfer characteristics of the solution space obtained from the two methods. These comparison results imply that the use of the diffusion model does not significantly affect the completeness of the solution space. This statement can be further supported by the comparison result in terms of the obtained minimum $\Delta v$ solution shown in Table \ref{tab4}, as the minimum $\Delta v$ solutions obtained from the proposed method are comparable to those obtained from the grid search method, for all three groups. These comparison results indicate that the proposed diffusion model-augmented method can generate a solution space with comparable transfer characteristics compared to the traditional grid search method. Then, on the basis of the comparable transfer characteristics, we further compare the two methods in terms of convergence rate and computational efficiency.

\begin{table}[h]
\centering
\caption{Comparison Result in terms of the Obtained Minimum $\Delta v$ Solution}\label{tab4}%
\begin{tabular}{@{}lll@{}}
\hline
Group & Grid Search, km/s & Proposed, km/s   \\
\hline
I    & 3.845075 & 3.845491 \\
II    & 3.692284 & 3.692340 \\
III & 3.644816 & 3.644891 \\
\hline
\end{tabular}
\end{table}

Table \ref{tab5} presents the convergence rate obtained from these two methods. It can be found that the use of the diffusion model-augmented method can further improve the convergence rate of the trajectory correction, indicating the high-quality initial guesses generated by the diffusion model. From Table \ref{tab5}, the relative improvement of the convergence rate over the traditional grid search method is 47.34-56.25$\%$ relatively, which can be considered as a remarkable advantage of the proposed method. Also, Table \ref{tab5} presents the wall-clock time of trajectory correction obtained from these two methods. It can be observed that the proposed method saves the wall-clock time more than 39.39$\%$ over the traditional method. These comparison results show that the proposed method takes advantage over the traditional grid search method in the convergence rate and computational efficiency using the high-quality initial guesses obtained from the diffusion model, further strengthening the motivation to combine the generative model with trajectory construction. To further explore the applicable range of the proposed method, we select three additional groups of departure/arrival orbital altitudes, and compare the diffusion model-augmented method with the traditional grid search method. The detailed results are presented in Appendix B, further illustrating that the preliminary applicable range of the proposed method with the samples generated from the diffusion model trained using the dataset in this paper ($h_i=167\text{ km}$ and $h_f=100\text{ km}$) is $h_i\in\left[167,\text{ }1500\right]\text{ }\text{ km}$ and $h_f\in\left[100,\text{ }1500\right]\text{ }\text{ km}$. A wider applicable range will be further explored in our future work. Moreover, the aforementioned phenomenon and construction results are presented in the Earth-Moon PCR3BP. When considering the Sun-Earth/Moon PBCR4BP or other higher-fidelity models, the distribution of $\left( {{\text{TOF}},{\text{ }}{\alpha _i}} \right)$ and transfer solution space might be more complex due to the additional construction parameters \cite{pinelli2023neural}. How to explore and utilize the solution space of the bi-impulsive Earth-Moon transfers in these models remains an open problem.
\begin{table}[h]
\centering
\caption{Comparison Result in terms of Convergence Rate}\label{tab5}%
\begin{tabular}{@{}lll@{}}
\hline
Group & Grid Search & Proposed   \\
\hline
I    & 1018001/2436480 (41.78$\%$) & 1906844/2920815 (65.28$\%$) \\
II    & 1132795/2436480 (46.49$\%$) & 2040593/2920815 (69.86$\%$) \\
III & 1212226/2436480 (49.75$\%$) & 2141073/2920815 (73.30$\%$) \\
\hline
\end{tabular}
\end{table}

\begin{table}[h]
\centering
\caption{Comparison Result in terms of Wall-Clock Time of Trajectory Correction}\label{tab6}%
\begin{tabular}{@{}lll@{}}
\hline
Group & Grid Search, s & Proposed, s   \\
\hline
I    & $2.6449\times 10^5$ & $1.6032\times 10^5$ \\
II    & $2.2349\times 10^5$ & $1.6083\times 10^5$ \\
III & $1.9142\times 10^5$ &  $1.1386\times 10^5$\\
\hline
\end{tabular}
\end{table}

\section{Conclusion}\label{sec5}
Construction of bi-impulsive Earth-Moon transfers with different departure/arrival orbital altitudes in the Earth-Moon planar circular restricted three-body problem is investigated in this paper. Firstly, the solution space of bi-impulsive transfers with different departure/arrival orbital altitudes is analyzed, in particular in terms of construction parameters, including departure phase angle at the Earth parking orbit, initial-to-circular velocity ratio, and time of flight. An interesting phenomenon about the time-of-flight distribution with respect to departure phase angle is discovered. Meanwhile, this type of distribution is similar between different departure/arrival orbital altitudes, which provides insight into the selection of the dataset to train a diffusion model. Then, a diffusion model is trained based on these phenomena, and a diffusion model-augmented grid search method is proposed to improve the convergence rate and computational efficiency of the traditional grid search method. For the construction of transfers with three settings of orbital altitudes, the proposed method improves the convergence rate by 47.34-56.25$\%$ and saves the wall-clock time by 39.39-40.52$\%$ over the traditional grid search method relatively, while ensuring comparable transfer characteristics. Applicable range of the proposed method is also preliminarily explored. This paper combines the training of a diffusion model with the multi-body dynamical phenomenon, augmenting the traditional grid search method and enabling convenient construction of Earth–Moon transfers with different departure/arrival orbital altitudes.

\section*{Declaration of competing interest}
The authors declare that they have no known competing financial interests or personal relationships that could have appeared to influence the work reported in this paper.

\section*{Acknowledgements}

This work was supported by the National Natural Science Foundation of China (Grant Nos. 12525204, 12372044, and U23B6002) and Beihang University Fund Program for Talent Development Support Scheme (Grant No. 4263023).

\appendix
\section*{Appendix A: Preliminary Interpretation for the Discontinuous Behavior of the TOF Distribution}

Taking Group I as an example, we provide a preliminary interpretation for the TOF discontinuous distribution. Firstly, we categorize the obtained transfer solutions (using the traditional grid search method) into transfers without the lunar gravity assist (LGA) and transfers with LGA. The criterion to identify a transfer with LGA is that there is more than one intersection with the disk $r_2=r_{\text{SOI}}$ in the transfer trajectory. The radius of the lunar sphere of influence (SOI) is set as $r_{\text{SOI}}=66243\text{ km}$ \cite{pan2022research}. The number of LGA can be calculated as $N_{\text{LGA}}=\left(N_{\text{Intersection}}-1\right)/2$, where $N_{\text{Intersection}}$ is the number of the intersections with the disk $r_2=R_{\text{SOI}}$ in the transfer. For the obtained solutions in Group I, the maximum $N_{\text{LGA}}$ is 4. The number of transfers with the corresponding $N_{\text{LGA}}$ is presented in Table \ref{num_LGA}.

\begin{table}[h]
\centering
\caption{The Number of Transfers with the Corresponding $N_{\text{LGA}}$}\label{num_LGA}%
\begin{tabular}{@{}ll@{}}
\hline
$N_{\text{LGA}}$ & Number of Transfers   \\
\hline
0    & 928156 \\
1    & 80256 \\
2    & 9264 \\
3    & 324 \\
4    & 1 \\
\hline
\end{tabular}
\end{table}

\begin{figure}
\centering
\includegraphics[width=0.99\textwidth]{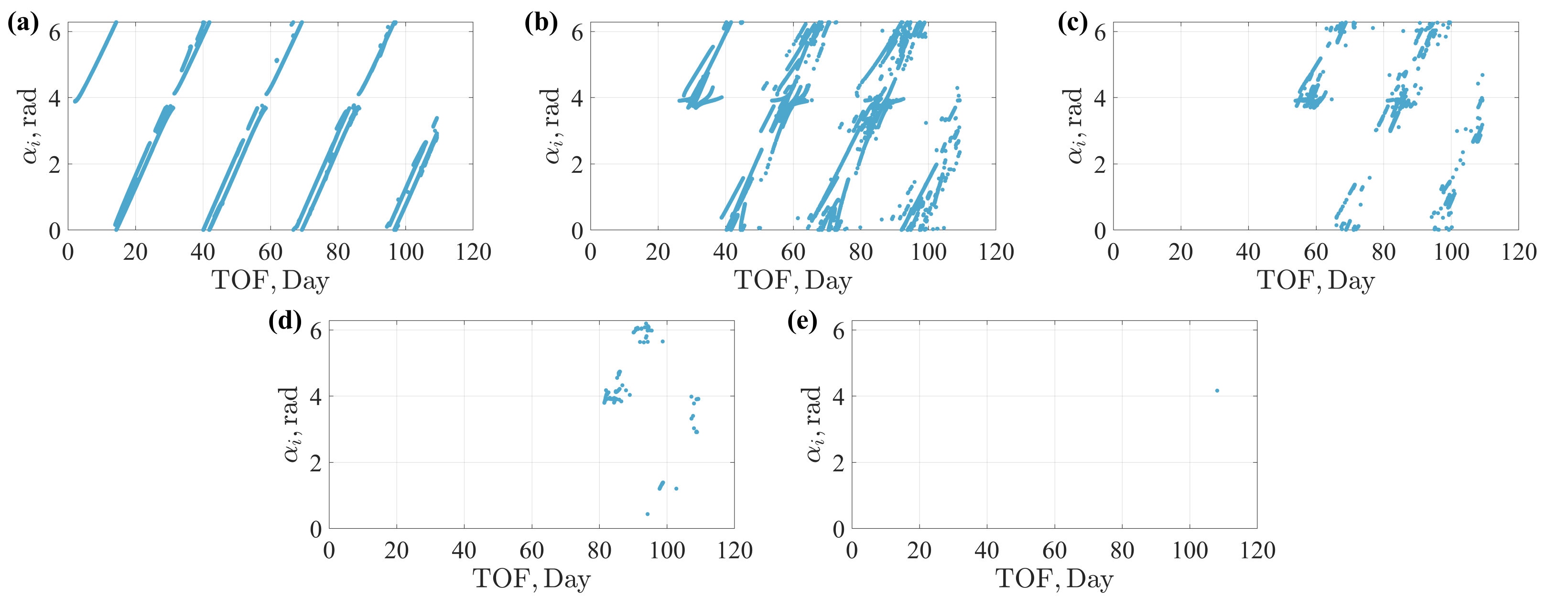}
\caption{The corresponding $\left(\text{TOF, }\alpha_i\right)$ maps for these five cases. (a) Transfers without LGA; (b) Transfers with LGA ($N_\text{LGA}=1$); (c) Transfers with LGA ($N_\text{LGA}=2$); (d) Transfers with LGA ($N_\text{LGA}=3$); (e) Transfers with LGA ($N_\text{LGA}=4$).}
\label{LGA}
\end{figure}

\begin{figure}
\centering
\includegraphics[width=0.4\textwidth]{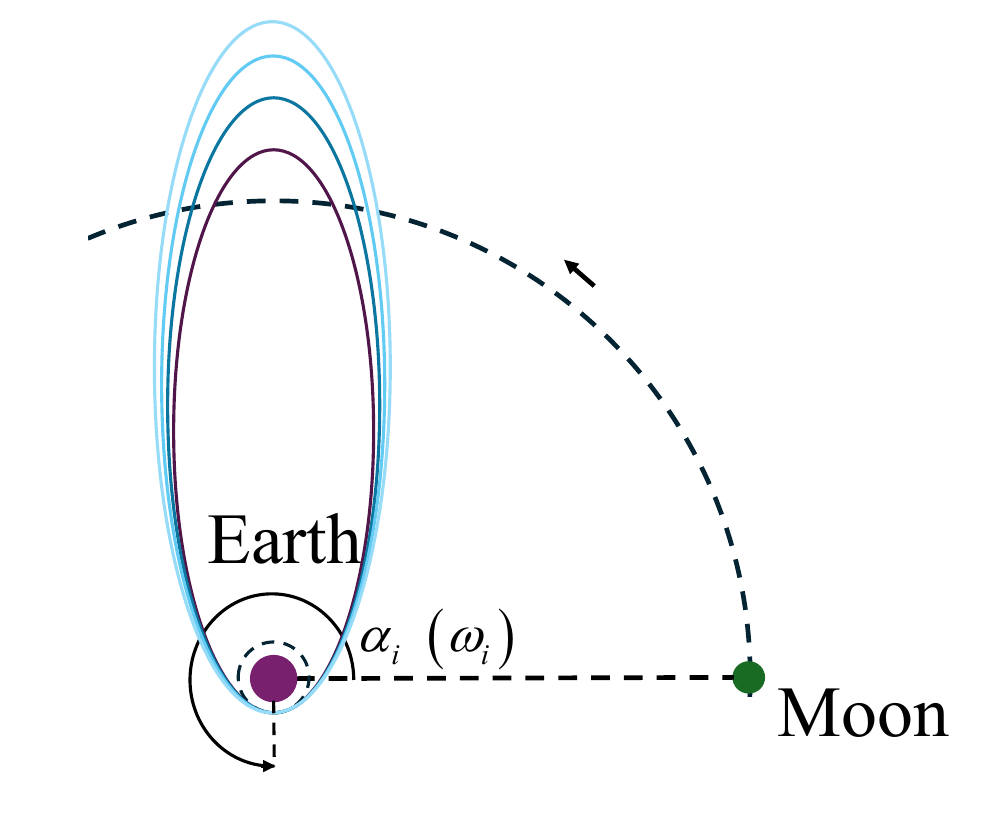}
\caption{A schematic to interpret the discontinuous behavior of the TOF distribution for transfers without LGA.}
\label{explanation}
\end{figure}

The corresponding $\left(\text{TOF, }\alpha_i\right)$ maps for these five cases are presented in Fig. \ref{LGA}. To analyze the overall $\left(\text{TOF, }\alpha_i\right)$ distribution, we focus on transfers without LGA and transfers with LGAs ($N_\text{LGA}=1\text{ and }2$). For transfers without LGA, before the trajectories enter the lunar SOI, the trajectories can be approximated as two-body orbits in the Earth-centered two-body problem. At the departure epoch of the bi-impulsive transfers, the phase angle of the Moon in the Earth-centered inertial (ECI) frame is set to 0. According to the assumption of the Earth-Moon PCR3BP mentioned in Section \ref{subsec2.1}, the Moon moves in a circular orbit around the Earth, and the corresponding period is $2\pi\frac{\text{TU}}{86400}$ Day (i.e., one month). We transform the transfers in the Earth-Moon rotating frame into those in the ECI frame. The transform equation is written as

\begin{equation}
\begin{gathered}
{\bm{X}_{in} = \bm{R}\left(\bm{X}-\bm{d}_{\text{Earth}}\right)} \hfill \\ 
{\bm{d}_{\text{Earth}} = \left[-\mu,\ 0,\ 0,\ 0\right]^{\text{T}}} \hfill \\ 
{\bm{R} =
\begin{bmatrix}
\bm{R}_{11} & \bm{O} \\ 
\bm{R}_{21} & \bm{R}_{11}
\end{bmatrix}} \hfill \\ 
{\bm{R}_{11} =
\begin{bmatrix}
\cos t & -\sin t \\ 
\sin t & \cos t
\end{bmatrix}} \hfill \\
{\bm{R}_{21} =
\begin{bmatrix}
-\sin t & -\cos t \\ 
\cos t & -\sin t
\end{bmatrix}} \hfill \\
\end{gathered}
\label{eq_ECI}
\end{equation}
where the subscript “\textit{in}” denotes the quantities in the ECI frame, and $t$ denotes the epoch during the transfer ($t\in \left[0,\text{ TOF}\right]$). Due to the tangential $\Delta v_i$, the trajectories (without LGA) before entering the lunar SOI can be approximated as high-eccentricity elliptic orbits ($\beta_i < \sqrt{2}$). Notably, the initial argument of perigee of these trajectories is $\alpha_i$. The argument of perigee $\omega$ of the trajectories can be calculated as

\begin{equation}
\begin{gathered}
 {{\bm{e}_{in}} = \frac{{\left( {{\bm{v}_{in}} \cdot {\bm{v}_{in}} - \frac{{1 - \mu }}{{{r_{in}}}}} \right){\bm{r}_{in}} - \left( {{\bm{r}_{in}} \cdot {\bm{v}_{in}}} \right){\bm{v}_{in}}}}{{1 - \mu }} } \hfill \\
  {\omega  = \bmod \left[ {\text{atan2}\left( {{e_{in,{\text{ }}y}},{\text{ }}{e_{in,{\text{ }}x}}} \right),{\text{ }}2\pi } \right]} \hfill \\
\end{gathered}
\label{eq_omega}
\end{equation}
where $\bm{X}_{in}=\left[\bm{r}_{in}^\text{T},\text{ }\bm{v}_{in}^\text{T}\right]^\text{T}$. Therefore, transfers with the same value of $\alpha_i$ have the same initial $\omega$. Therefore, for approximately two-body orbits, there is a relatively concentrated time interval each month (i.e., the period of the Moon around the Earth) when approaching the Moon, creating opportunities for Earth-Moon transfers, which provides an interpretation of this phenomenon for the transfers without LGA. A schematic of this geometric relationship is shown in Fig. \ref{explanation}.

\begin{figure}
\centering
\includegraphics[width=0.4\textwidth]{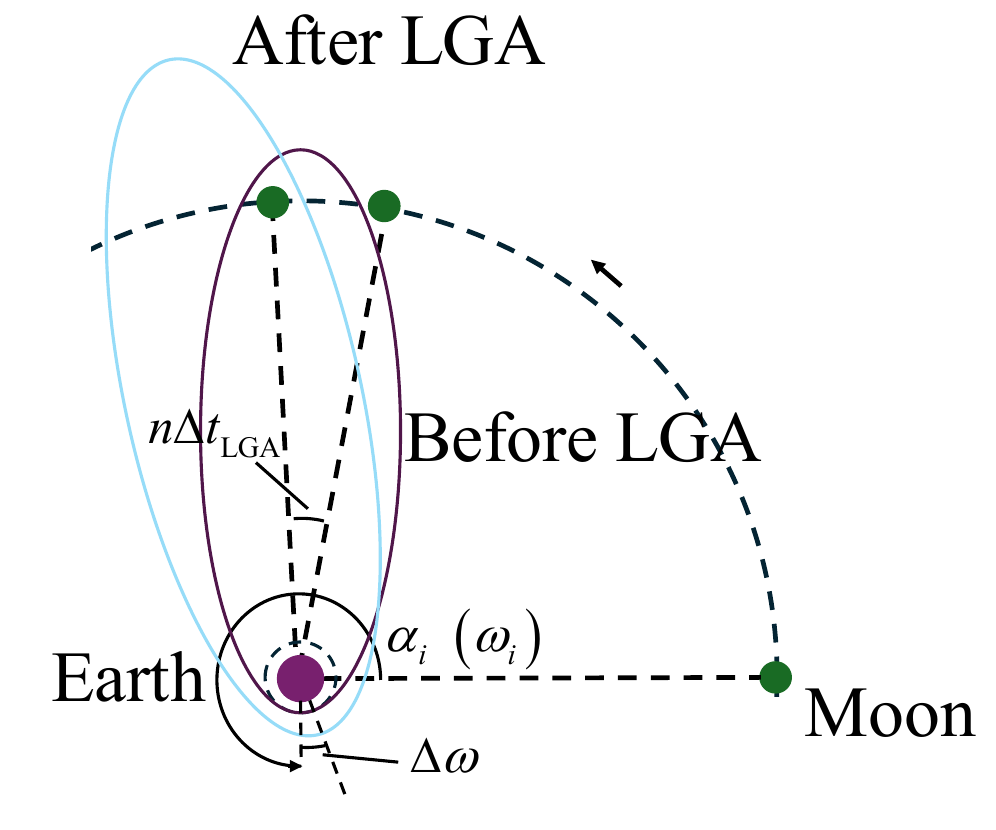}
\caption{A schematic for transfers with LGA.}
\label{LGA_LGA}
\end{figure}

\begin{figure}
\centering
\includegraphics[width=0.3\textwidth]{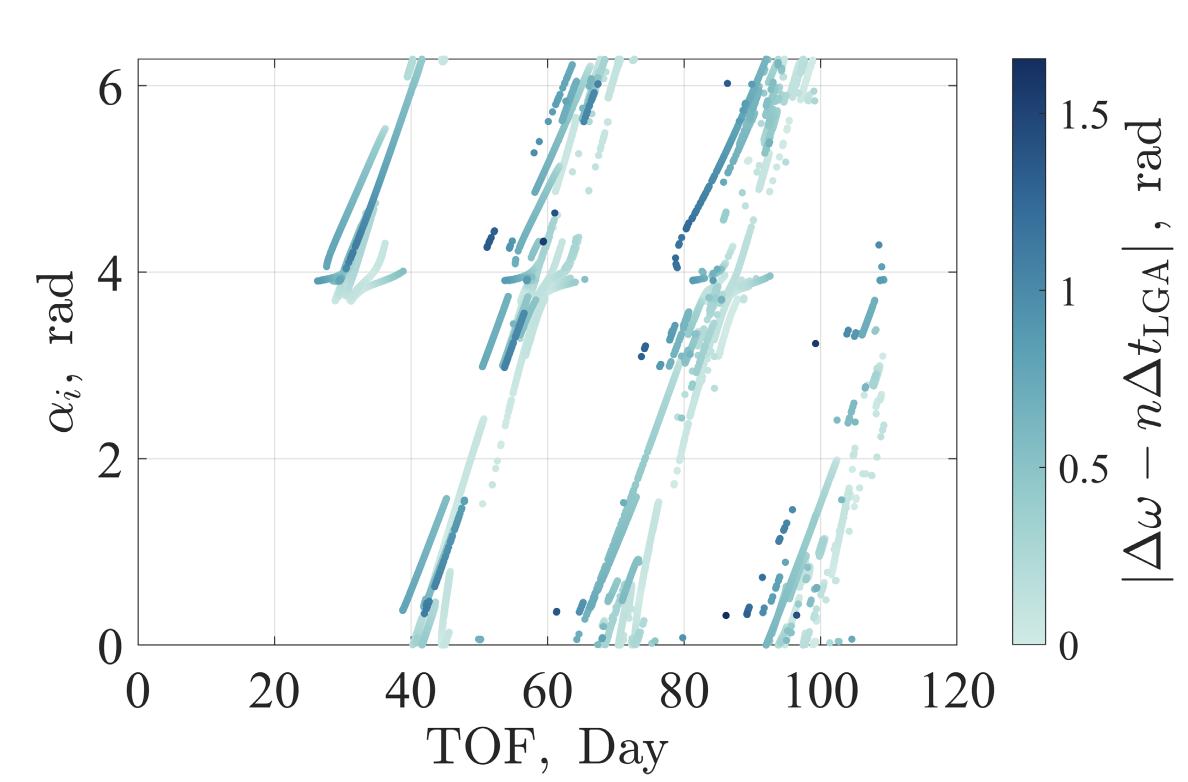}
\caption{The distribution of $\left|\Delta \omega-n\Delta t_\text{LGA}\right|$ in terms of $\left(\text{TOF, }\alpha_i\right)$ ($N_\text{LGA}=1$).}
\label{error_1}
\end{figure}

\begin{figure}
\centering
\includegraphics[width=0.6\textwidth]{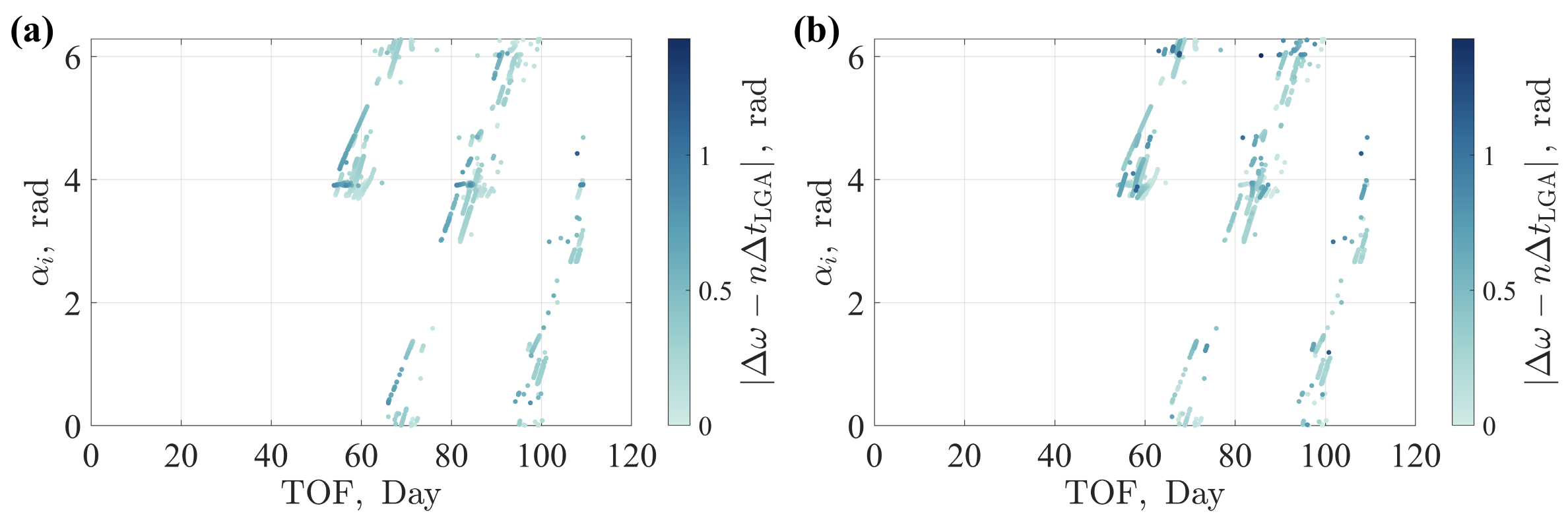}
\caption{The distribution of $\left|\Delta \omega-n\Delta t_\text{LGA}\right|$ in terms of $\left(\text{TOF, }\alpha_i\right)$ ($N_\text{LGA}=2$). (a) First LGA; (b) Second LGA.}
\label{error_2}
\end{figure}

For transfers with LGA, we focus on the transfers corresponding to $N_\text{LGA}=1$ and $N_\text{LGA}=2$. Due to the LGA, $\omega$ can have a relatively pronounced variation. Then, we use the numerical method to determine the effects of this variation on the discontinuous behavior of the TOF distribution. For transfers with one LGA, we record the time of LGA $\Delta t_\text{LGA}$ (approximated as the time interval during the trajectory within the lunar SOI before finally arriving in the circular Moon target orbit) and the variation of the argument of perigee $\Delta \omega$ after the LGA (the corresponding values of $\omega$ are calculated by the states entering the lunar SOI and departing the lunar SOI). An LGA represents one opportunity to approach the Moon. After the LGA, to maintain the discontinuous behavior of the TOF distribution, the value of $\Delta \omega-n\Delta t_\text{LGA}$ (here $\Delta t_\text{LGA}$ adopts the dimensionless time unit) should be relatively low to maintain the relative geometric relationship between the Moon and the approximate two-body orbit (after the LGA), as shown in Fig. \ref{LGA_LGA}. The statement is also suitable for transfers with two LGAs. Then, we report the computation results supporting this statement. Table \ref{epsilon} presents the median and P95 of $\left|\Delta \omega-n\Delta t_\text{LGA}\right|$ for each LGA for transfers with LGA, and the corresponding distribution in terms of $\left(\text{TOF, }\alpha_i\right)$ is presented in Figs. \ref{error_1} and \ref{error_2}. From Table \ref{epsilon} and these figures, it can be observed that most of the values of $\Delta \omega-n\Delta t_\text{LGA}$ are relatively low. Therefore, the discontinuous behavior of the TOF distribution also exists for transfers with LGA. Moreover, it is $\Delta \omega-n\Delta t_\text{LGA}$ that possibly causes the banded distribution of $\alpha_i$ with respect to TOF. The aforementioned discussion provides a preliminary interpretation for the discontinuous behavior of the TOF distribution from the perspective of the geometric relationship between the transfer trajectories and the Moon in the ECI frame. However, the specific mechanisms of the LGA still require further exploration. It will be the focus of our future work.

\begin{table}[h]
\centering
\caption{Median and P95 of $\left|\Delta \omega-n\Delta t_\text{LGA}\right|$ for Each LGA}\label{epsilon}%
\begin{tabular}{@{}llll@{}}
\hline
$N_{\text{LGA}}$ & First/Second LGA & Median of $\left|\Delta \omega-n\Delta t_\text{LGA}\right|$, rad & P95 of $\left|\Delta \omega-n\Delta t_\text{LGA}\right|$, rad  \\
\hline
1 & -- & 0.163306 & 0.678437 \\
2 & First LGA & 0.196563 & 0.611277 \\
2 & Second LGA & 0.248229 & 0.633783 \\
\hline
\end{tabular}
\end{table}

\section*{Appendix B: Preliminary Exploration for the Applicable Range of the Proposed Method}

\begin{figure}
\centering
\includegraphics[width=0.75\textwidth]{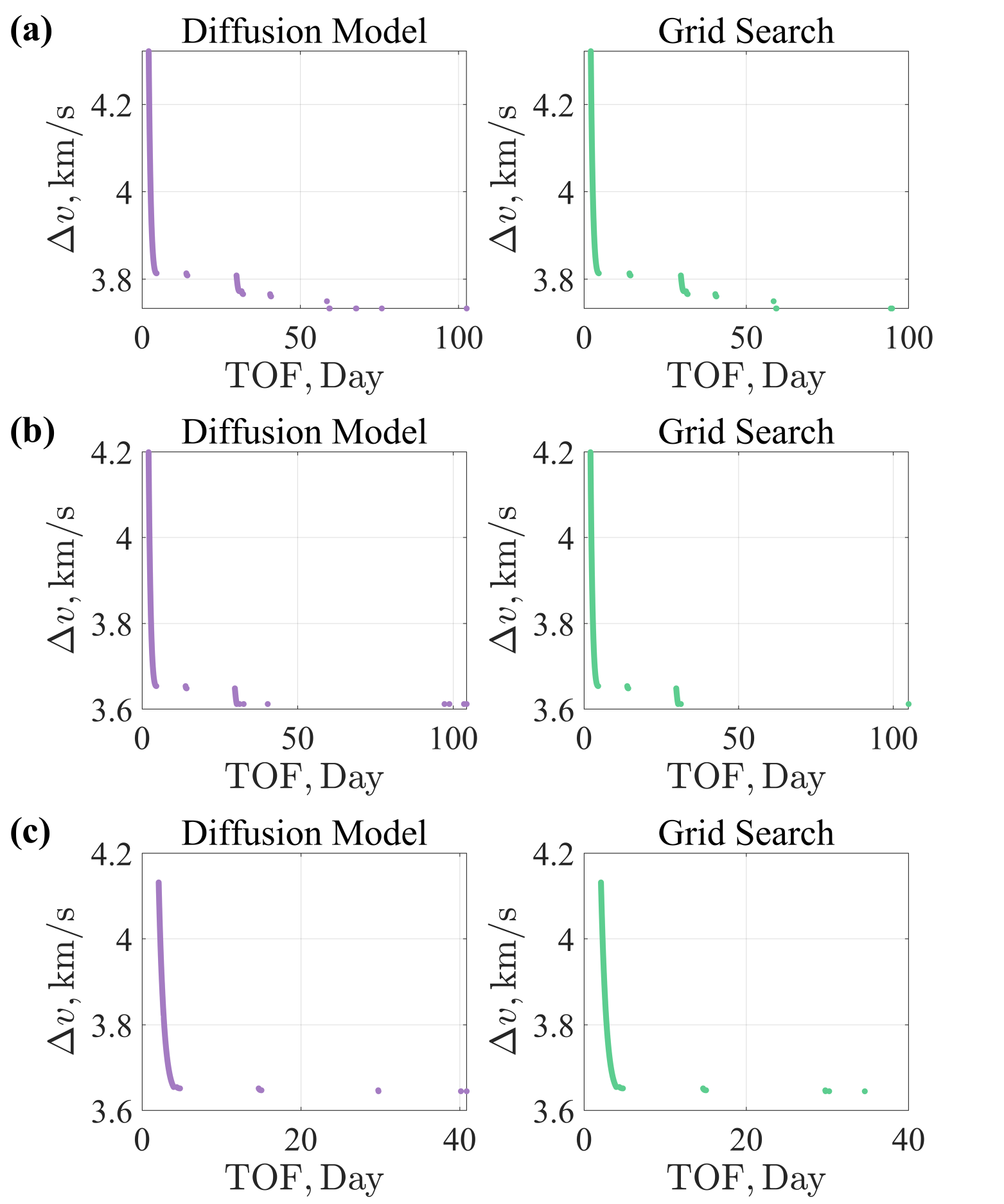}
\caption{Pareto fronts obtained from the two methods. (a) Group IV; (b) Group V; (c) Group VI.}
\label{comparison_pareto_group_IV_VI}
\end{figure}

Here, we select three additional groups of departure/arrival orbital altitudes and further explore the applicable range of the proposed method. These three groups are set as

\begin{equation} \begin{gathered}
  {{\text{Group IV:     }}{h_i} = 500{\text{ km,     }}{h_f} = 500{\text{ km}}} \hfill \\
 {{\text{Group V:    }}{h_i} = 1000{\text{ km,     }}{h_f} = 1000{\text{ km}}} \hfill \\
  {{\text{Group III:   }}{h_i} = 1500{\text{ km,   }}{h_f} = 1500{\text{ km}}} \hfill \\
\end{gathered} 
\label{eq6_new}
\end{equation}
We further perform the diffusion model-augmented method (with the generated data in Section \ref{subsubsec3.2.4}) and the traditional grid search method to construct the corresponding bi-impulsive Earth-Moon transfers. The comparison in terms of the Pareto fronts of the obtained solutions are presented in Fig. \ref{comparison_pareto_group_IV_VI}. The obtained minimum $\Delta v$ solutions for these three groups are presented in Table \ref{tab_dv_group_4_6}. The results shown in Fig. \ref{comparison_pareto_group_IV_VI} and Table \ref{tab_dv_group_4_6} also suggest comparable transfer characteristics for these three groups obtained from the proposed method and the traditional grid search method. On the basis of the comparable transfer characteristics, we further investigate the convergence rate and computational efficiency, which are presented in Tables \ref{tab5_new} and \ref{tab6_new}. From these two tables, it can be found that for these three groups, the proposed method improves the convergence rate and computational efficiency compared to the traditional grid search method, suggesting that the proposed method is also applicable to these three groups. Based on the aforementioned discussion, here, we present a preliminary applicable range of the proposed diffusion model-augmented construction method (with the dataset selected in Section \ref{subsec3.1}): $h_i\in\left[167,\text{ }1500\right]\text{ }\text{ km}$ and $h_f\in\left[100,\text{ }1500\right]\text{ }\text{ km}$. A wider applicable range can be obtained from further exploration by investigating more groups of departure/arrival orbital altitudes in future work.

\begin{table}[h]
\centering
\caption{Comparison Result in terms of the Obtained Minimum $\Delta v$ Solution (Groups IV-VI)}\label{tab_dv_group_4_6}%
\begin{tabular}{@{}lll@{}}
\hline
Group & Grid Search, km/s & Proposed, km/s   \\
\hline
IV    & 3.732600 & 3.732604 \\
V    & 3.611911 & 3.611938 \\
VI & 3.514419 & 3.514416 \\
\hline
\end{tabular}
\end{table}

\begin{table}[h]
\centering
\caption{Comparison Result in terms of Convergence Rate (Groups IV-VI)}\label{tab5_new}%
\begin{tabular}{@{}lll@{}}
\hline
Group & Grid Search & Proposed   \\
\hline
IV    & 1141289/2436480 (46.84$\%$) & 2043005/2920815 (69.95$\%$) \\
V    & 1355696/2436480 (55.64$\%$) & 2202407/2920815 (75.40$\%$) \\
VI & 1339787/2436480 (54.99$\%$) & 2288337/2920815 (78.35$\%$) \\
\hline
\end{tabular}
\end{table}

\begin{table}[h]
\centering
\caption{Comparison Result in terms of Wall-Clock Time of Trajectory Correction (Groups IV-VI)}\label{tab6_new}%
\begin{tabular}{@{}lll@{}}
\hline
Group & Grid Search, s & Proposed, s   \\
\hline
IV    & $2.7073\times 10^5$ & $1.8234\times 10^5$ \\
V    & $2.2718\times 10^5$ & $1.4376\times 10^5$ \\
VI & $1.6253\times 10^5$ &  $1.1601\times 10^5$\\
\hline
\end{tabular}
\end{table}



\printcredits

\end{document}